\documentclass[12pt]{article} 

\usepackage[english]{babel}
\usepackage[a4paper,top=3cm,bottom=2cm,left=2cm,right=2cm,marginparwidth=1.75cm]{geometry}

\usepackage[style = numeric,
            sorting = nty]{biblatex}

\usepackage{float}
\usepackage{wrapfig}
\usepackage{amsmath}
\usepackage{mathtools}
\usepackage{authblk}
\usepackage{enumitem}
\usepackage{wrapfig}

\newcommand{\at}[2][]{#1|_{#2}}

\usepackage{xcolor}
\newcommand{\rv}[1]{\textcolor{black}{#1}}

\title{Compact implicit high resolution numerical method for solving transport problems with sorption isotherms
\thanks{This work was supported by the grants VV-MVP-24-0116, VEGA 1/0314/23 and APVV-23-0186.}}

\author[1]{Dagmar \v{Z}\'AKOV\'A}
\author[2]{Peter FROLKOVI\v{C}}
\affil[1]{Department of Mathematics and Descriptive Geometry, Slovak University of Technology in Bratislava (dagmar.zakova@stuba.sk)}
\affil[2]{Department of Mathematics and Descriptive Geometry, Slovak University of Technology in Bratislava (peter.frolkovic@stuba.sk)}

\begin{document}
\date{}
\maketitle

\textbf{Abstract:} 
This study investigates numerical methods to solve nonlinear transport problems characterized by various sorption isotherms with a focus on the Freundlich type of isotherms. We describe and compare second order accurate numerical schemes, focusing on implicit methods, to effectively model transport phenomena without stability restriction on the choice of time steps. Furthermore, a high resolution form of the method is proposed that limits a priori the second order accurate scheme towards first order accuracy to keep the values of numerical solutions in a physically acceptable range.

Through numerical experiments, we demonstrate the effectiveness of high resolution methods in minimizing oscillations near discontinuities, thereby enhancing solution plausibility. The observed convergence rates confirm that the second order accurate schemes achieve expected accuracy for smooth solutions and that they yield significant improvements when compared with the results of the first order scheme. As the computational cost of the compact implicit method seems to be comparable to similar explicit ones with a clear profit of unconditional stability, this research provides a practical tool toward numerical simulations of nonlinear transport phenomena applicable in various fields such as contaminant transport in porous media or column liquid chromatography.


\section{Introduction}
\label{sec-introduction}

When solving numerically time dependent partial differential equations (PDEs), one has to decide what kind of time discretization is appropriate for each particular mathematical model and its applications. The two basic classes of time discretization techniques are the so-called explicit ones and the implicit ones. The simplest time discretization of the time derivative in the PDE can be obtained, e.g., by forward finite difference (the explicit Euler method) which typically leads to an explicit definition of the values for the numerical solution in a new time level depending only on approximations evaluated in the old (known) time level. This is a big advantage and so, if possible, appropriate higher order explicit time discretization methods, typically of a Runge-Kutta type, are used for the time discretization.

It is well known that such explicit time discretization techniques require some restriction on the choice of time steps. Nevertheless, if this restriction on the size of time steps due to stability is in accordance with the requirement on its size due to accuracy demand, the stability restriction should not cause practical problems.

The implicit time discretizations do not offer, in general, any explicit definitions of the numerical solution, and to obtain it one has to solve a system of algebraic equations. Although the property of no stability restriction on the time step is an attractive feature for numerical simulators, the cost of algebraic solvers can be too high to compensate for it. 

The main difference between explicit and implicit time discretization methods is significantly relaxed in the case when the term with the time derivative in the PDE is nonlinear. In such cases also explicit time discretization methods can lead to algebraic equations for the values of numerical solution, although, in general, of a simpler form than for the implicit time discretization. 

In this paper, we are interested in numerical solution of nonlinear transport equations in a conservative form where the time derivative in the model is applied to a nonlinear function of the solution. One can easily show that such problems are formally equivalent to nonlinear hyperbolic PDEs \cite{frolkovic_semi-analytical_2006,javeedEfficientAccurateNumerical2011} even in the case of systems of PDEs \cite{burger_linearly_2018, donatImplicitExplicitWENO2018}, but such a transformation does not need to be available in an analytical form, therefore, the origin formulation of PDEs has to be solved. 

For situations involving simplified initial and boundary conditions, analytical solutions are available to some extent \cite{frolkovic_semi-analytical_2006,qamarAnalyticalSolutionsMoment2014, kaczmarskiAnalyticalNumericalSolutions2024}. In general, one has to use numerical methods to solve such models when explicit time discretization is preferred for the advection part of the model \cite{kacurSolutionContaminantTransport2005,javeedEfficientAccurateNumerical2011,donatImplicitExplicitWENO2018,burger_linearly_2018,zafarDiscontinuousGalerkinScheme2021,donat_weno_2024}.

In \cite{frolkovic2023high,zakova_numerical_2024}, we have successfully applied high resolution compact implicit numerical methods for the solution of standard hyperbolic PDEs. The main advantage of such methods is that they require nonlinear algebraic solvers only for small systems opposite to standard implicit methods, where large algebraic systems of coupled nonlinear equations should be solved.

In this work, we apply for the first time such a numerical scheme to the transport equation with nonlinear sorption isotherms. Such models are used to describe the transport of contaminants in porous media when contaminant adsorption on pore skeletons must be taken into account \cite{frolkovicNumericalSimulationContaminant2016}. Other interesting applications of such models are in liquid chromatography technology \cite{guMathematicalModelingScaleUp2015}. We show that to solve related representative mathematical models, the computational cost of the compact implicit scheme is comparable to explicit schemes, but with the clear advantage of no stability requirement for the implicit scheme.

The paper is organized as follows. In Section 2, we introduce the mathematical models. The numerical schemes are described in Section 3.  The two-dimensional case is detailed in Section 4. Finally, extensive numerical experiments are presented in Section 5.

\section{Mathematical model}
\label{sec-mathmodel}

We consider the representative nonlinear transport problem in the form 
\begin{equation}
    \label{eq_general}
    \partial_{t} F(u) + \partial_{x} u= 0,
\end{equation}
with $u = u(x, t) \geq 0$ being the unknown function for $t\in(0,T)$ and $x\in(x_L,x_R)\subset R$. The initial condition is defined by 
\begin{equation}
    \label{initialcondition}
    u(x, 0)=u^0(x) \,.
\end{equation}
As $F'(u)>0$ for $u \in R$, the Dirichlet boundary condition is defined only at the left boundary,
\begin{equation}
    \label{boundarycondition}
    u(x_L,t) = u_{x_L} (t) \,,
\end{equation}
and no boundary condition is required at the right boundary.

The mathematical model can represent the transportation of contaminant \cite{brusseau1995effect,kacurSolutionContaminantTransport2005, frolkovic_semi-analytical_2006,frolkovicNumericalSimulationContaminant2016} or liquids in chromatography \cite{guMathematicalModelingScaleUp2015,schmidt-traubPreparativeChromatography2020} with the sorption isotherm 
\begin{equation}
    \label{psi}
    \Psi(u) = F(u) - u \,.
\end{equation}
One of the most frequently encountered type of (nonlinear) sorption isotherms is the Freundlich type,
\begin{equation}
    \label{Freundlich}
    \Psi(u) = au^p\,, \quad  a, p>0 \,.
\end{equation}
In this paper, we will focus on a single choice $a=1$.

The solution to (\ref{eq_general}) discussed in this paper is based on the theoretical results \cite{leveque_finite_2004}, which ensure the existence and uniqueness of so-called entropy solutions for hyperbolic problems. Such theoretical framework has been established for hyperbolic equations in the form
\begin{equation}
    \label{eq_withinvers}
    \partial_{t} q + \partial_{x} f(q)= 0 \,,  
\end{equation}
that can be obtained with the transformation 
\begin{equation}
    \label{eq_withinvers_2}
    q = F(u)\,, \,\,  u=f(q) \rv{\,,}
\end{equation} 
where $f=F^{-1}$ is the inverse function of F.

In general, the function $f$ may not be available analytically. 
In the case where $f$ is known, the problem can be readily addressed using established methods for hyperbolic problems \cite{leveque_finite_2004}.

The equation \eqref{eq_general} may be extended using the velocity $v(x)$, known for each value of $x$, into the form
\begin{equation}
    \label{eq_general_withvel}
    \partial_{t} F(u) +  \partial_{x} (v(x)u)= 0 \,,
\end{equation}
where the velocity may take on both positive and negative values \rv{and the Dirichlet boundary condition will be extended also to the right boundary
\begin{equation}
    \label{boundarycondition_right}
    u(x_R,t) = u_{x_R} (t) \,.
\end{equation}}

\section{Numerical schemes}
\label{sec-numscheme}

The numerical discretization to solve the equation \eqref{eq_general_withvel} is done in space and time using the notation $t^n=n\tau$, $n=0,1,...N$ for a chosen $N$ and $\tau>0$, and the spatial discretization is based on the finite volume method \cite{leveque_finite_2004}. The computational domain is divided into finite volumes with a regular grid in the form \rv{ $I_{i} = (x_{i-1/2},x_{i+1/2}) $, where $x_{i-1/2} = x_L +(i-1)h$ with $i=1, 2, ..., M$ }for the chosen $M$ and $h=(x_R-x_L)/M$.

To obtain the numerical scheme, we use the finite volume method, with the numerical approximation of space and time averaged values  $\bar q_{i}^n$ and \rv{ $\bar u_{i+1/2}^{n+1/2}$} as in, e.g., \cite{shu_essentially_1998,leveque_finite_2004},
\begin{equation}
    \label{average}
        \begin{split}
        \bar q_i^n  
        & \coloneqq \frac{1}{h} \int_{I_i} F\big(u(x,t^n)\big) \,dx  \,,\\
        \bar u_{i+1/2}^{n+1/2} &\coloneqq  \frac{1}{\tau} \int_{t^n}^{t^{n+1}} u(x_{i+1/2},t) \,dt \,.
        \end{split}
\end{equation}

To solve \eqref{eq_general_withvel}, we need to treat the positive and negative velocities separately. Let
\begin{equation}
    \label{velocitysplit}
    v= v^+ + v^-, \hspace{1cm}  \text{with  } \, v^+ \coloneqq \max(0,v), \,\, v^{-} \coloneqq \min (0,v) \,,
\end{equation}
with the values of $v$ defined at the edges of finite volume $v_{i+1/2} = v(x_{i+1/2})$. Note that the equation \eqref{eq_general} is obtained from \eqref{eq_general_withvel} simply by setting $v \equiv 1$.

After integrating the differential equation (\ref{eq_general_withvel}) over $I_i \times (t^n,t^{n+1})$ and using the definitions (\ref{average}), we obtain the scheme 
\begin{equation}
    \label{exactscheme}
    \bar q_i^{n+1} - \bar q_i^n + \frac{\tau}{h} (v_{i+1/2}\bar u_{i+1/2}^{n+1/2} - v_{i-1/2}\bar u_{i-1/2}^{n+1/2}) = 0 \,.
\end{equation}
The numerical solution of \eqref{exactscheme} is obtained by approximating the values $\bar u_{i+1/2}^{n+1/2}$ using $U_{i+1/2}^{n+1/2} \approx \bar u_{i+1/2}^{n+1/2}$ \cite{shu_essentially_1998,leveque_finite_2004} and $Q_{i}^n \approx \bar q_{i}^n 
$ denote an approximation of the averaged values $\bar q_i^{n+1}$ in (\ref{average}). Consequently, the numerical scheme takes the form
\begin{equation}
    \label{numericalschemeGeneral}
    Q_i^{n+1} - Q_i^n + \frac{\tau}{h} (v_{i+1/2}U_{i+1/2}^{n+1/2} - v_{i-1/2}U_{i-1/2}^{n+1/2}) = 0 \,.
\end{equation}
To complete the description of the numerical scheme, we have to define the values $U_{i+1/2}^{n+1/2}$ using the numerical values $U_i^{n}$ and $Q_i^{n+1}=F(U_i^{n+1})$.

\subsection{First order schemes}
\label{sec-1storder}

For convenience and a comparison, we give notation for the first order accurate explicit and implicit numerical schemes. To obtain the first order accurate explicit numerical scheme, we use the upwind approach together with the splitting \eqref{velocitysplit}, and the scheme \cite{leveque_finite_2004} will become
\begin{align}
    \label{firstorder_Explicit_withvel_scheme}
    Q_i^{n+1} = Q_i^n - \frac{\tau}{h} \Big[ &v_{i+1/2}^+ U_{i}^{n}   -  v_{i-1/2}^+ U_{i-1}^{n}    \nonumber \\ 
                      +  & v_{i+1/2}^- U_{i+1}^{n} - v_{i-1/2}^- U_{i}^{n}   \Big]  \,,
\end{align}
\rv{where the explicit structure of the scheme is clearly evident.}

It is important to note that when employing an explicit numerical scheme \eqref{firstorder_Explicit_withvel_scheme}, the challenge of solving the nonlinear problem does not vanish and is still present in terms $Q_i^{n+1}=F(U_i^{n+1})$. The nonlinearity $F(u)=u+u^p$ still requires resolution, which necessitates the use of an iterative method to obtain the solution, such as Newthon's method, together with the necessity to adhere to the stability restrictions imposed by explicit numerical schemes.

In case of using the implicit numerical scheme, one gets \rv{
\begin{align}
    \label{firstorder_Implicit_withvel_scheme}
    Q_i^{n+1} - Q_i^n + \frac{\tau}{h} \Big[ &v_{i+1/2}^+ U_{i}^{n+1}   -  v_{i-1/2}^+ U_{i-1}^{n+1}    \nonumber \\ 
                      +  & v_{i+1/2}^- U_{i+1}^{n+1} - v_{i-1/2}^- U_{i}^{n+1}   \Big] = 0 \,,
\end{align}
where the nonlinearity is readily apparent in the previously mentioned relation $Q_i^{n+1} = F(U_i^{n+1})$.}

The discretization scheme creates a system of algebraic equations \eqref{firstorder_Implicit_withvel_scheme} that can be solved iteratively, using the fast sweeping method \cite{zhao2005fast}, where each iteration is given by a nonlinear Gauss-Seidel iteration.
The approach of the fast sweeping method indicates that Gauss-Seidel iteration is applied sequentially in alternating index directions across the computational domain. Specifically, we $"sweep"$ from each end of the domain as follows:
\begin{align}
    \label{sweepiter1}
    \text{First sweep:} \quad &i = 1, \dots, M-1, \\
    \label{sweepiter2}
    \text{Second sweep:} \quad &i = M-1, \dots, 1. 
\end{align}
At each of these $"sweeps"$, we solve scalar nonlinear algebraic equation (\ref{firstorder_Implicit_withvel_scheme}) for each index $i$ using the Newton's method. Generally, to reduce the absolute error of the solution, more than two fast sweeping iterations may be required for each computational time step. Note that in the case of $v(x)>0$, only the first sweep is necessary (\ref{sweepiter1}). Conversely, when $v(x)<0$, the second sweep alone is sufficient (\ref{sweepiter2}).

\subsection{Second order accurate scheme}
\label{2ndorder}

In this section, we will focus on the derivation of the second order accurate compact implicit numerical scheme and its high resolution formulation.

To obtain the second order of accuracy \cite{duraisamy_implicit_2007} of (\ref{numericalschemeGeneral}), considering for a moment only a positive velocity $v(x)>0$, when we approximate the values $U_{i+1/2}^{n+1/2}$ applying a finite Taylor series,
\begin{align}
    \label{taylor}
    u(x_{i+1/2},t^{n+1/2}) = u(x_{i},t^{n+1}) + \frac{h}{2} \partial_x u(x_{i}, t^{n+1})  - \frac{\tau}{2}  \partial_t u(x_{i}, t^{n+1}) 
     + \mathcal{O}(h^2, \tau^{2}) \,,
\end{align}   
where the numerical values would then be defined as
\begin{align}
    \label{desired_numerical_values}
        U_{i+1/2}^{n+1/2} &= U_{i}^{n+1} + \frac{h}{2} \partial_x U_{i}^{n+1} - \frac{\tau}{2} \partial_t U_{i}^{n+1}  \, ,
\end{align}
with the approximations $\partial_x U_{i}^{n+1}$ and $\partial_t U_{i}^{n+1}$ expressed as a linear combination of two different approximation choices \cite{frolkovic2023high} using parameter $\boldsymbol{\omega} \in [0,1]$ \rv{
\begin{align}
        \partial_x U_{i}^{n+1} & = \frac{1}{h} \big( \omega (U_{i}^{n+1} - U_{i-1}^{n+1}) + (1-\omega) (U_{i+1}^{n+1} - U_{i}^{n+1}) \big) \,,\nonumber \\
        \label{derivatives_approx}
        \partial_t U_{i}^{n+1} & = \frac{1}{\tau} \big( \omega (U_{i}^{n+1} - U_{i}^{n}) + (1-\omega) (U_{i+1}^{n+1} - U_{i+1}^{n}) \big) \,.
\end{align}
Consequently, combining \eqref{desired_numerical_values} with \eqref{derivatives_approx}, we obtain the compact implicit numerical approximation
\begin{equation}
    \label{U_approx_2ndord_implicit}
        U_{i+1/2}^{n+1/2} = U_{i}^{n+1} - \frac{1}{2}\big(\omega(U_{i-1}^{n+1} - U_{i}^{n}) + (1-\omega)(U_{i}^{n+1} - U_{i+1}^{n})\big) \,.
\end{equation}}


To deal with the oscillations that may arise due to discontinuities in the solution, we propose a high resolution form of the scheme. Such scheme incorporates a variable parameter $\boldsymbol{\omega}$, and a limiter $\boldsymbol{l}$ to control numerical oscillations and to ensure stability \cite{frolkovic2023high, zakova_numerical_2024}. By adapting $\boldsymbol{\omega}$ and employing the limiter $\boldsymbol{l}$, the scheme effectively balances accuracy and non-oscillatory behavior near discontinuities. Consequently, the proposed method is transformed into the following form,
\begin{align}
    \label{scheme_2nd_implicit}
        Q_i^{n+1} - Q_i^n + \frac{\tau}{h} \Big[ v_{i+1/2} \Big(&U_{i}^{n+1} - \frac{l_i}{2}\big(\omega_{i}(U_{i-1}^{n+1} - U_{i}^{n}) + (1-\omega_{i})(U_{i}^{n+1} - U_{i+1}^{n})\big) \Big) \nonumber \\
        - v_{i-1/2} \Big(&U_{i-1}^{n+1} - \frac{l_{i-1}}{2}\big(\omega_{i-1}(U_{i-2}^{n+1} - U_{i-1}^{n}) + (1-\omega_{i-1})(U_{i-1}^{n+1} - U_{i}^{n})\big)\Big)  \Big] = 0 \,.
\end{align}

Again, it is necessary to handle the nonlinearity present in \rv{$Q$, since $Q=F(U)$}, and to determine the values of $\boldsymbol{\omega}$ and $\boldsymbol{l}$. To do that, we use the WENO (Weighted Essentially Non-Oscillatory) approximations, together with the predictor-corrector algorithm, which was described in detail in \cite{zakova_numerical_2024, frolkovic2023high}.

In case of solving the equation (\ref{eq_general_withvel}) with the velocity that changes the sign, we need to extend the values $\boldsymbol{\omega} = (\boldsymbol{\omega}^+, \boldsymbol{\omega}^-)$ and $\boldsymbol{l}=(\boldsymbol{l}^+, \boldsymbol{l}^-)$ and the scheme becomes high resolution
\begin{align}
    \label{scheme_2nd_implicit_withvel}
    Q_i^{n+1} - Q_i^n +     
       \frac{\tau}{h}\Big[ v_{i+1/2}^+\big(&U_{i}^{n+1}   - \frac{l_i^+}{2}\big(\omega_{i}^+   (U_{i-1}^{n+1} - U_{i}^{n}  ) + (1-\omega^+_i)    (U_{i}^{n+1}   - U_{i+1}^{n}) \big) \big) \nonumber \\
     -                            v_{i-1/2}^+\big(&U_{i-1}^{n+1} - \frac{l_{i-1}^+}{2}\big(\omega_{i-1}^+ (U_{i-2}^{n+1} - U_{i-1}^{n}) + (1-\omega_{i-1}^+)(U_{i-1}^{n+1} - U_{i}^{n}  ) \big) \big) \nonumber \\
     + v_{i+1/2}^-\big(&U_{i+1}^{n+1} - \frac{l_{i+1}^-}{2}\big(\omega_{i+1}^- (U_{i+2}^{n+1} - U_{i+1}^{n}) + (1-\omega_{i+1}^-)(U_{i+1}^{n+1} - U_{i}^{n}  ) \big) \big) \nonumber \\
     -                            v_{i-1/2}^-\big(&U_{i}^{n+1}   - \frac{l_i^-}{2}\big(\omega_{i}^  - (U_{i+1}^{n+1} - U_{i}^{n}  ) + (1-\omega_{x}^-  )(U_{i}^{n+1}   - U_{i-1}^{n}) \big)\big) \Big]  \, = 0 \,,
\end{align}
which, again, may be solved using the fast sweeping method described in Section \ref{sec-1storder}.

\section{Two dimensional case}
\label{sec-2D}

We consider the model for the unknown function $u=u(x,y,t)$, where $(x,y) \in R^2$ with $x\in(x_L,x_R)\subset R$, $y\in(y_L,y_R)\subset R$, in the form 
\begin{equation}
    \label{2D_withvel}
    \partial_{t} F(u) +  \partial_{x} (v(x,y) u) +  \partial_{y} (w(x,y) u) = 0,
\end{equation}
where $\Vec{v}= \Vec{v}(x,y)=(v(x,y),w(x,y))$ represents a divergence free velocity field, which is a known function, and the initial condition is given  $u(x,y,t)=u^0(x,y)$ and the Dirichlet boundary conditions, if prescribed, are denoted by
\begin{align}
    \label{boundary2D}
    u(x_L,y,t) = u_{x_L} (y,t)\,, &\quad u(x,y_L,t) = u_{y_L} (x,t) \,, \nonumber \\
    u(x_R,y,t) = u_{x_R} (y,t)\,, &\quad u(x,y_R,t) = u_{y_R} (x,t) \,.
\end{align}

\rv{
For simplicity of the notation, we assume a squared computational domain ($x,y\in(x_L,x_R)$) that is divided into finite volumes of the form $I_{i,j} = (x_{i-1/2},x_{i+1/2}) \times (y_{j-1/2},y_{j+1/2})$ with the regular square grid, where $x_{i-1/2} = x_L + (i-1)h$ and $y_{j-1/2} = x_L + (j-1)h$ with $i,j=1, 2, ..., M$ for the chosen $M$ and $h=(x_R-x_L)/M$.}

The process is similar to the one described in Section \ref{sec-mathmodel}, with the addition of one spatial variable. The detailed derivation is also described in \cite{AACEE_2023_Zakova, zakova_numerical_2024} up to second order accuracy. 

The final numerical scheme for the velocity splitting 
\begin{align}
    \label{2D_velocitysplit}
    v=v^+ + v^-, \hspace{1cm} & \text{with  } \, v^+ \coloneqq \max(0,v), \,\, v^{-} \coloneqq \min (0,v)\,, \nonumber \\
    w=w^+ + w^-, \hspace{1cm} & \text{with  } \, w^+ \coloneqq \max(0,w), \,\, w^{-} \coloneqq \min (0,w) \rv{\,,}
\end{align}
will look like 
\begin{align}
    \label{2D_numericalscheme}
         Q_{i,j}^{n+1} - Q_{i,j}^n  +\frac{\tau}{h} \Big[&v^+_{i+1/2,j}{U}_{i+1/2,j}^{n+1/2} - v^+_{i-1/2,j}{U}_{i-1/2,j}^{n+1/2} + w^+_{i,j+1/2}{U}_{i,j+1/2}^{n+1/2} - w^+_{i,j-1/2}{U}_{i,j-1/2}^{n+1/2} \nonumber\\
         + &v^-_{i+1/2,j}{U}_{i+1/2,j}^{n+1/2} - v^-_{i-1/2,j}{U}_{i-1/2,j}^{n+1/2}  + w^-_{i,j+1/2}{U}_{i,j+1/2}^{n+1/2} - w^-_{i,j-1/2}{U}_{i,j-1/2}^{n+1/2} \Big] = 0 \rv{\,,}
\end{align}
where the numerical values $U_{i+1/2,j}^{n+1/2}$, $U_{i,j+1/2}^{n+1/2}$, and $U_{i-1/2,j}^{n+1/2}$ $U_{i,j-1/2}^{n+1/2}$ need to be determined differently for positive and negative velocities, based on the upwind approach, similarly to Section \ref{sec-numscheme}. 

The simplest first order accurate implicit numerical scheme takes the form
\begin{align}
    \label{scheme1ord}
        Q_{i,j}^{n+1} - Q_{i,j}^n +
         \frac{\tau}{h} \Big[ &v^+_{i+1/2,j}{U}_{i,j}^{n+1} - v^+_{i-1/2,j}{U}_{i-1,j}^{n+1} +  w^+_{i,j+1/2}{U}_{i,j}^{n+1} - w^+_{i,j-1/2}{U}_{i,j-1}^{n+1} \nonumber \\
         + &v^-_{i+1/2,j}{U}_{i+1,j}^{n+1} - v^-_{i-1/2,j}{U}_{i,j}^{n+1} + w^-_{i,j+1/2}{U}_{i,j+1}^{n+1} - w^-_{i,j-1/2}{U}_{i,j}^{n+1} \Big] = 0 \,,
\end{align}
and can be, again, solved iteratively using the fast sweeping method described in Section \ref{sec-mathmodel} with the nonlinearity also in the $Q_{ij}^{n+1}$ term. This time we use four different $"sweeps"$ from each corner of a rectangular domain where each fast sweeping iteration is given by one Gauss-Seidel iteration applied sequentially: 
\begin{align}
    \label{2D_sweeps}
    \text{First sweep:} \quad i = 1, ..., M \ &, \ \ j = 1, ..., M \nonumber \\
    \text{Second sweep:} \quad i = M, ..., 1 \ &, \ \ j = 1, ..., M \nonumber \\
    \text{Third sweep:} \quad i = M, ..., 1 \ &, \ \ j = M, ..., 1 \nonumber \\ 
    \text{Fourth sweep:} \quad i = 1, ..., M \ &, \ \ j = M, ..., 1  
\end{align}
where one has to solve a scalar nonlinear algebraic equation for each pair of indices $i$ and $j$ in the sweeps.

\rv{
The high resolution compact implicit numerical scheme will then use the approximations for the values ${U}_{i+1/2,j}^{n+1/2}$ and ${U}_{i,j+1/2}^{n+1/2}$ in \eqref{2D_numericalscheme} as
\begin{align}
         {U}_{i+1/2,j}^{n+1/2} & = U_{i,j}^{n+1} - \frac{l_{i,j}^{x,+}}{2}\big(\omega_{i,j}^{x,+} (U_{i-1,j}^{n+1} -  U_{i,j}^{n}) + (1-\omega_{i,j}^{x,+})(U_{i,j}^{n+1} - U_{i+1,j}^{n})\big) \,, \nonumber \\
    \label{2D_secondorderscheme_positive}
         {U}_{i,j+1/2}^{n+1/2} & = U_{i,j}^{n+1} - \frac{l_{i,j}^{y,+}}{2}\big(\omega_{i,j}^{y,+} (U_{i,j-1}^{n+1} - U_{i,j}^{n}) + (1-\omega_{i,j}^{y,+})(U_{i,j}^{n+1} - U_{i,j+1}^{n})\big) \,, 
\end{align}
which stand next to the positive velocities $v^+_{i+1/2,j}$ and $w^+_{i,j+1/2}$; and the approximations ${U}_{i+1/2,j}^{n+1/2}$ and ${U}_{i,j+1/2}^{n+1/2}$ for negative cases 
\begin{align}
        {U}_{i+1/2,j}^{n+1/2} & = U_{i+1,j}^{n+1} - \frac{l_{i+1,j}^{x,-}}{2}\big(\omega_{i+1,j}^{x,-} (U_{i+2,j}^{n+1} - U_{i+1,j}^{n}) + (1-\omega_{i+1,j}^{x,-})(U_{i+1,j}^{n+1} - U_{i,j}^{n}) \big) \,, \nonumber \\
    \label{2D_secondorderscheme_negative}
        {U}_{i,j+1/2}^{n+1/2} & = U_{i,j+1}^{n+1} - \frac{l_{i,j+1}^{y,-}}{2}\big(\omega_{i,j+1}^{y,-} (U_{i,j+2}^{n+1} - 
           U_{i,j+1}^{n}) + (1-\omega_{i,j+1}^{y,-})(U_{i,j+1}^{n+1} - U_{i,j}^{n}) \big) \,,  
\end{align}
which stand next to the velocities $v^-_{i+1/2,j}$ and $w^-_{i,j+1/2}$, }with the parameter $\boldsymbol{l}=(\boldsymbol{l}^{x,\pm}, \boldsymbol{l}^{y,\pm})  \in [0,1]$, and $\boldsymbol{\omega} = (\boldsymbol{\omega}^{x,\pm}, \boldsymbol{\omega}^{y,\pm})  \in [0,1]$, different for each $I_{i,j}$.  For a detailed derivation of the scheme and definitions of parameters, we refer to \cite{AACEE_2023_Zakova, zakova_numerical_2024}.

\section{Numerical experiments}
\label{sec-experiments}

In this part, we will discuss the experiments using the Freundlich type of isotherm (\ref{Freundlich}) $q = F(u) = u + u^p$, for different choices of intervals 
\begin{equation}
    \label{p_choices}
        p\begin{cases}
        p \in (0,1) &  \\
        p \in (1,2) &  \\
        p \geq 2
        \end{cases}
\end{equation}

To demonstrate the order of accuracy, we compute the discrete $L_1$ norm of the error $E$, calculated as 
\begin{equation}
    \label{error}
    E = h \sum \limits_{i=1}^M  \lvert U_{i}^N - \bar u_{i}^N \lvert \,
\end{equation}
and the experimental order of convergence (EOC) of the numerical results for the examples with the known exact solution (Section \ref{subsec_1D_withouvel}). As a reference for comparing the results, the maximal Courant number ($C_{max}$) is calculated, taking into account the derivative of $F$ as the retardation factor
\begin{equation}
    \label{F_derivative}
    \frac{\partial F(u)}{\partial u} \at[\biggl]{u} = 1 + pu^{p-1} \geq 1 \,, \quad u \ge 0 \,.
\end{equation}
Considering that the retardation factor is always greater or equal to $1$, the maximum Courant number will be defined simply by
\begin{equation}
    \label{Courantnumber}
    C_{max} =  \frac{\tau}{h} \max \limits_{i} \Bigl\{ \max \bigl\{ |v(x_{i-1/2})|, |v(x_{i+1/2})| \bigr\}  \Bigr\} \,.
\end{equation}

As mentioned, Newton's method was proposed as an iterative method to solve the nonlinear problem. This may lead to difficulty when the derivative of $F(u)$ must be computed for the method, especially in cases where $u = 0$. For such cases, regularization of $u$ is proposed in the following manner:
\begin{equation}
    u_i = \max \bigl\{ u_i, 10^{-6} \bigr\} \,.
\end{equation}

The numerical methods and graphical outputs in this work were implemented using the Python programming language \cite{CS-R9526}. 

\subsection{One dimensional experiment}
\label{subsec_1D_withouvel}

\begin{figure}[h!]
    \begin{center}
    \includegraphics[width=17cm]{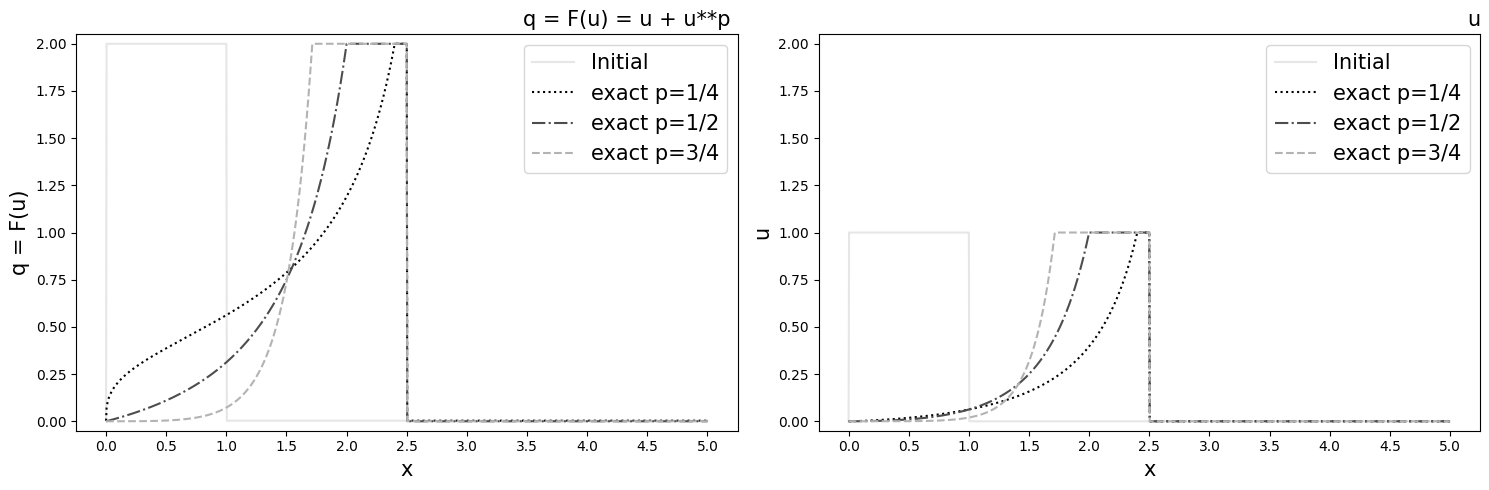}
    \caption{The exact solution of $q$ and $u$ for $p=1/4,1/2,3/4$ for the final time $T=3$ with the initial condition (\ref{InitialCondition_discontinuous}).}
    \label{fig:Exact_pmensieako1}
    \end{center}
\end{figure}

\begin{figure}[h!]
    \begin{center}
    \includegraphics[width=17cm]{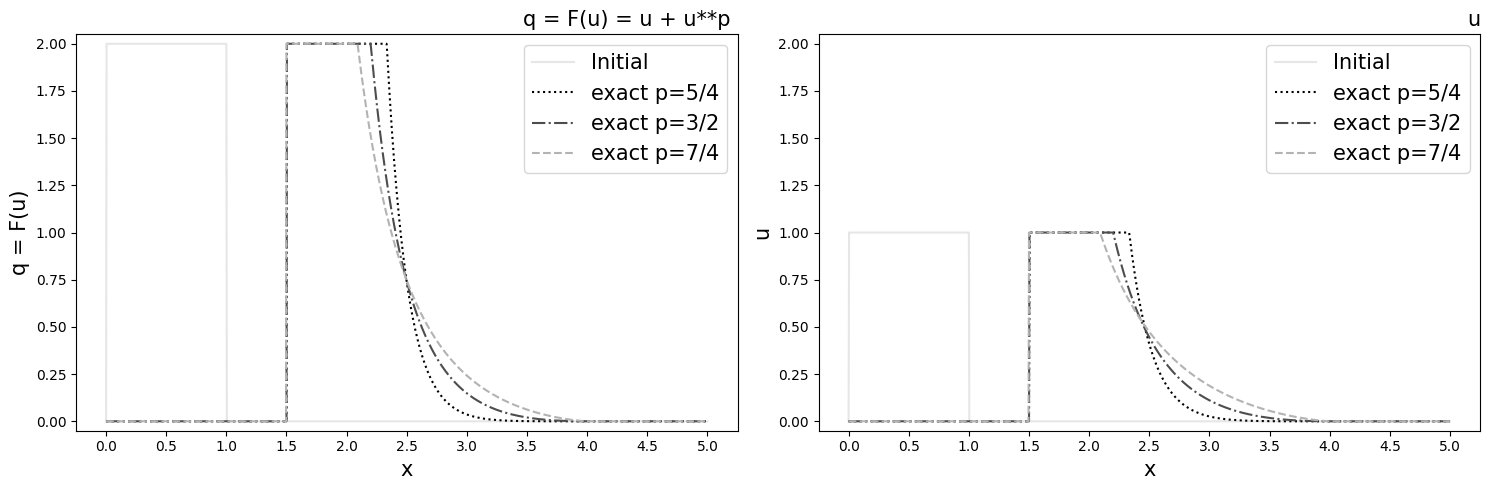}
    \caption{The exact solution of $q$ and $u$ for $p=5/4,3/2,7/4$ for the final time $T=3$ with the initial condition (\ref{InitialCondition_discontinuous}).}
    \label{fig:Exact_pvacsieako1}
    \end{center}
\end{figure}

\begin{figure}[h!]
    \begin{center}
    \includegraphics[width=17cm]{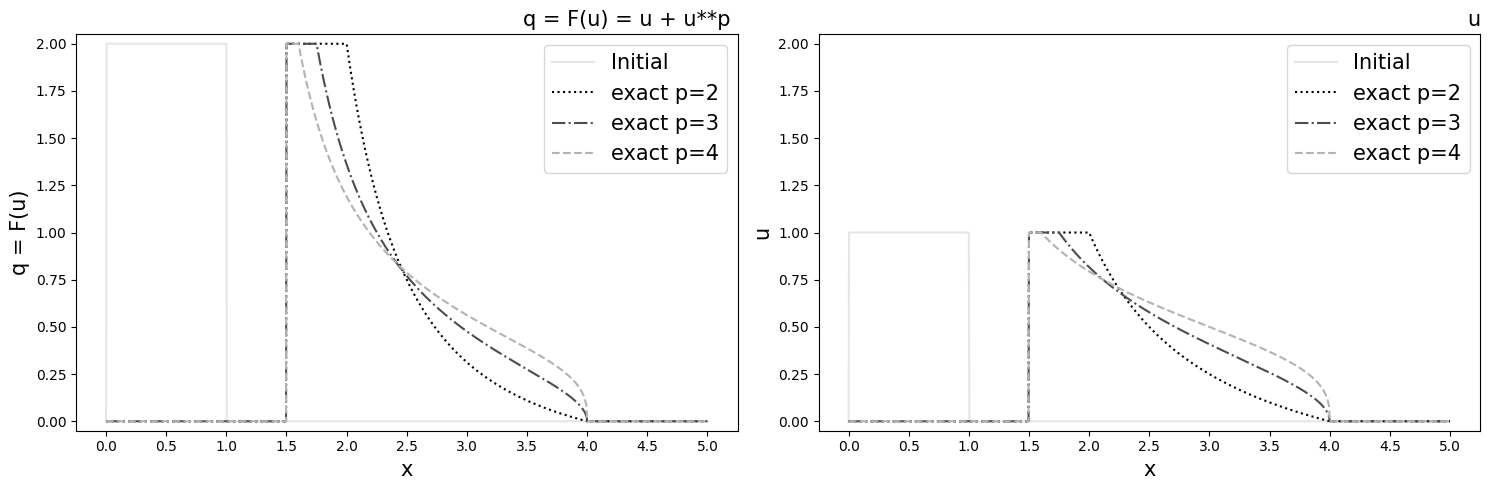}
    \caption{The exact solution of $q$ and $u$ for $p=2,3,4$ for the final time $T=3$ with the initial condition (\ref{InitialCondition_discontinuous}).}
    \label{fig:Exact_pvacsieako2}
    \end{center}
\end{figure}

For the first example we choose the velocity $v(x)=1$ and the initial condition defined for $x\in[0,5]$ as a discontinuous function in a form
\begin{equation}
    \label{InitialCondition_discontinuous}
    u^0(x) =
        \begin{cases}
        1 & \text{for } 0 < x < 1 \\
        0 & \text{otherwise}
        \end{cases}
\end{equation}
with the boundary condition $u_{0}(t)=0$ and we will compute the solution for the final time $T=3$. The exact solutions are described in \cite{frolkovic_semi-analytical_2006} for each of the choices of $p$ (\ref{p_choices}) and are shown in Figures \ref{fig:Exact_pmensieako1} - \ref{fig:Exact_pvacsieako2}. 

Firstly, we use the computational domain $x\in[0.5,1.5]$ and the time period $t\in[2,3]$, where the solution is free of discontinuities that could otherwise affect the accuracy order. To use the full stencil of the scheme in this example for every inner grid point, especially for the second order accurate schemes, we set exact values not only at the boundary point, but also in a neighboring point outside of the computational interval.

In Table \ref{TAB:Results_cn-0.25}, we also present the CPU time to illustrate the computational duration differences between first order and second order accurate schemes, both implicit and explicit ones. The second order accurate explicit numerical scheme used for comparison is described in detail in \cite{qiu_finite_2003}. Applying the second order accurate numerical scheme to the problem \eqref{eq_general_withvel} still presents the necessity of solving the nonlinear problem due to $q = F(u)$.

For each scheme, the nonlinear algebraic equation must be solved at every time step and for each finite volume. The comparison of numerical solutions is conducted using $N=M*2$ and $N=M$ time steps, for $M=320, 640, 1280, 2560$ finite volumes, with the Courant numbers being $C_{max} = 0.25$ and $C_{max} = 0.5$, respectively, to ensure consistency in the results. The value $p=1/2$ in (\ref{Freundlich}) was chosen for this example. 

\begin{table}[h!]
    \begin{center}
    \begin{tabular}{c c | c c c | c c c }
    \hline
    \multicolumn{2}{c}{} & \multicolumn{6}{c}{First order, $C_{max} = 0.25$}  \\
    \multicolumn{2}{c}{} & \multicolumn{3}{c}{Explicit} & \multicolumn{3}{c}{Implicit}  \\
    \hline
    $M $ & $N$  & E                   & EOC  & CPU time  & E                    & EOC  & CPU time  \\ 
    \hline
    320  & 640  & 8.01 $\cdot 10^{-4}$ & -    & 2.54     & 1.05 $\cdot 10^{-3}$ & -    & 2.62      \\
    640  & 1280 & 4.02 $\cdot 10^{-4}$ & 0.99 & 11.00    & 5.28 $\cdot 10^{-4}$ & 0.99 & 12.24     \\
    1280 & 2560 & 2.01 $\cdot 10^{-4}$ & 0.99 & 42.20    & 2.64 $\cdot 10^{-4}$ & 0.99 & 44.36     \\
    2560 & 5120 & 1.00 $\cdot 10^{-4}$ & 0.99 & 161.15   & 1.32 $\cdot 10^{-4}$ & 0.99 & 169.18    \\
    \hline
    \multicolumn{2}{c}{} &  \multicolumn{6}{c}{Second order, $C_{max} = 0.25$}  \\
    \multicolumn{2}{c}{} & \multicolumn{3}{c}{Explicit} & \multicolumn{3}{c}{Implicit, $\omega=1/2$}  \\
    \hline
    $M $ & $N$  & E                    & EOC  & CPU time & E                  & EOC  & CPU time  \\ 
    \hline
    320  & 640  & 4.85 $\cdot 10^{-6}$ & -    & 5.42   & 2.94 $\cdot 10^{-6}$ & -    & 6.32      \\
    640  & 1280 & 1.23 $\cdot 10^{-6}$ & 1.98 & 24.00  & 7.58 $\cdot 10^{-7}$ & 1.95 & 21.25     \\
    1280 & 2560 & 3.09 $\cdot 10^{-7}$ & 1.99 & 95.15  & 1.92 $\cdot 10^{-7}$ & 1.97 & 85.97     \\
    2560 & 5120 & 7.76 $\cdot 10^{-8}$ & 1.99 & 360.77 & 4.84 $\cdot 10^{-8}$ & 1.98 & 333.50    \\
    \hline 
    \hline\\
    \multicolumn{2}{c}{} & \multicolumn{6}{c}{First order, $C_{max} = 0.5$}  \\
    \multicolumn{2}{c}{} & \multicolumn{3}{c}{Explicit} & \multicolumn{3}{c}{Implicit}  \\
    \hline
    $M $ & $N$  & E                   & EOC  & CPU time  & E                    & EOC  & CPU time  \\ 
    \hline
    320  & 320  & 6.75 $\cdot 10^{-4}$ & -    & 1.34     & 1.17 $\cdot 10^{-3}$ & -    & 1.40      \\
    640  & 640 & 3.38 $\cdot 10^{-4}$ & 0.99 & 5.023    & 5.91 $\cdot 10^{-4}$ & 0.99 & 7.02     \\
    1280 & 1280 & 1.69 $\cdot 10^{-4}$ & 0.99 & 21.85    & 2.96 $\cdot 10^{-4}$ & 0.99 & 23.70     \\
    2560 & 2560 & 8.47 $\cdot 10^{-5}$ & 0.99 & 84.70    & 1.48 $\cdot 10^{-4}$ & 0.99 & 90.62    \\
    \hline
    \multicolumn{2}{c}{} &  \multicolumn{6}{c}{Second order, $C_{max} = 0.5$}  \\
    \multicolumn{2}{c}{} & \multicolumn{3}{c}{Explicit} & \multicolumn{3}{c}{Implicit, $\omega=1/2$}  \\
    \hline
    $M $ & $N$  & E                    & EOC  & CPU time & E                  & EOC  & CPU time  \\ 
    \hline
    320  & 320  & 7.01 $\cdot 10^{-6}$ & -    & 3.17   & 2.67 $\cdot 10^{-6}$ & -    & 2.56     \\
    640  & 640 & 1.76 $\cdot 10^{-6}$ & 1.98 & 12.34  & 6.93 $\cdot 10^{-7}$ & 1.95 & 11.52     \\
    1280 & 1280 & 4.42 $\cdot 10^{-7}$ & 1.99 & 48.87  & 1.76 $\cdot 10^{-7}$ & 1.97 & 43.50    \\
    2560 & 2560 & 1.10 $\cdot 10^{-7}$ & 1.99 & 193.64 & 4.44 $\cdot 10^{-8}$ & 1.98 & 175.11    \\
    \hline 
    \end{tabular}
    \end{center}
    \caption{A comparison of the numerical errors, the EOCs, and the CPU times for the first and second order accurate schemes, both explicit and implicit, for a case with a smooth initial condition, using different maximum Courant numbers, $C_{max} = 0.25, 0.5$, with $p = 1/2$.}
    \label{TAB:Results_cn-0.25}
\end{table}

Table \ref{TAB:Results_cn-0.25} clearly demonstrates that the use of an explicit numerical scheme is not particularly advantageous, as the computational time and error levels are nearly identical to those of the implicit numerical scheme. This similarity arises because even when employing an explicit scheme, for example the first order accurate (\ref{firstorder_Explicit_withvel_scheme}), the nonlinear problem must still be solved at each step, negating the typical computational simplicity associated with explicit methods. 

\begin{figure}[h!]
    \begin{center}
    \includegraphics[width=17cm]{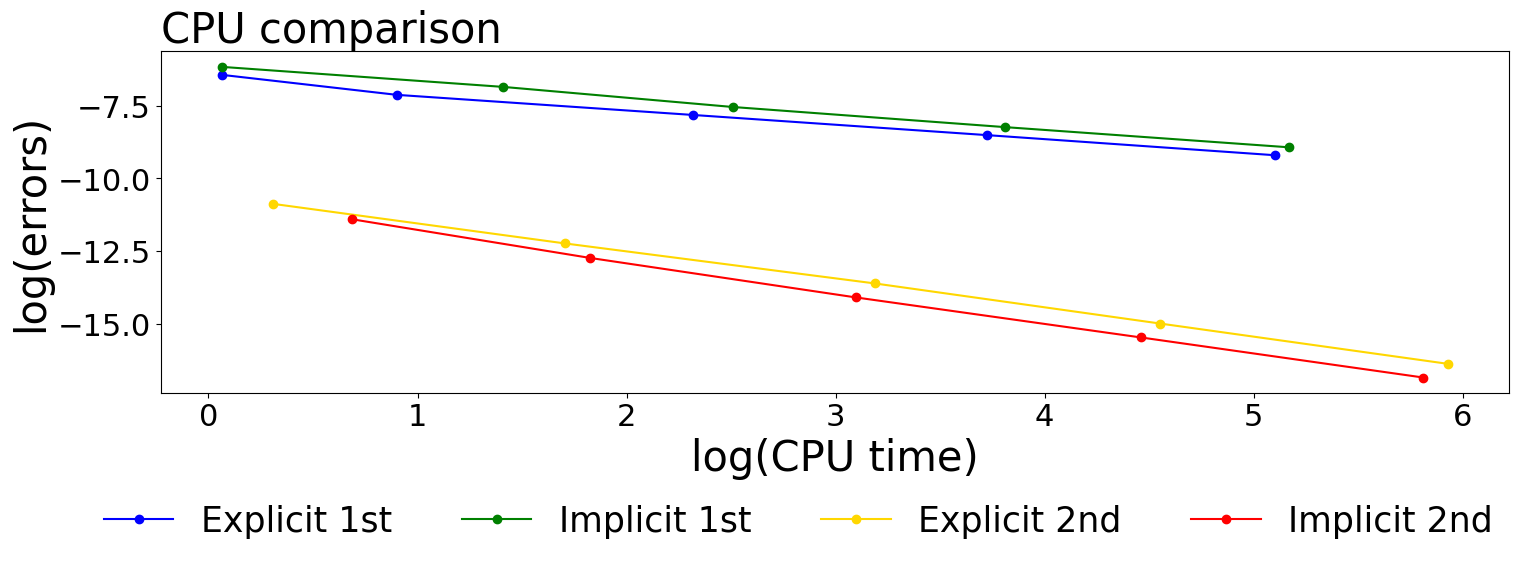}
    \caption{CPU time versus errors in log scale for the first and second order accurate numerical scheme (explicit and implicit) with $C_{max}=0.25$.}
    \label{fig:cpu_cn0.25}
    \end{center}
\end{figure}

\begin{figure}[h!]
    \begin{center}
    \includegraphics[width=17cm]{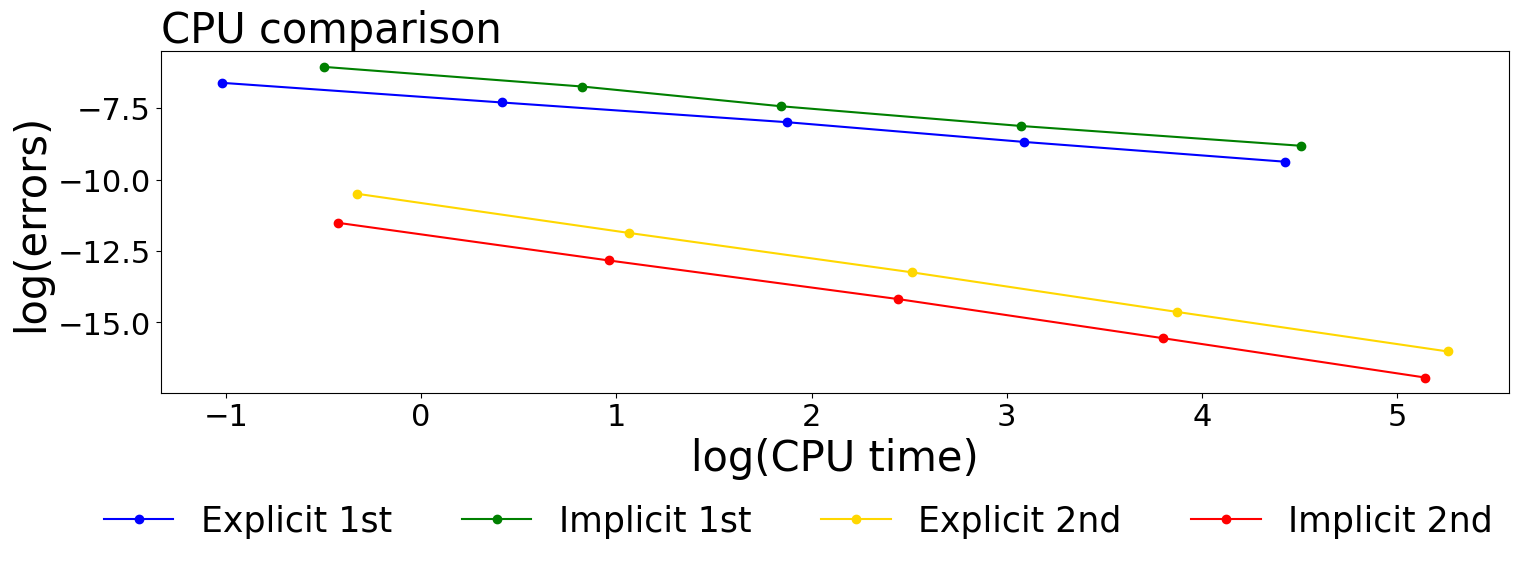}
    \caption{CPU time versus errors in log scale for the first and second order accurate numerical scheme (explicit and implicit) with $C_{max}=0.5$.}
    \label{fig:cpu_cn0.5}
    \end{center}
\end{figure}

\rv{ 
This can also be observed in Figures \ref{fig:cpu_cn0.25} for $C_{max} = 0.25$ and in Figure \ref{fig:cpu_cn0.5} for $C_{max} = 0.5$, which present the comparison of CPU time versus errors on a logarithmic scale.
}

\begin{figure}[h!]
    \begin{center}
    \includegraphics[width=17cm, height = 5cm]{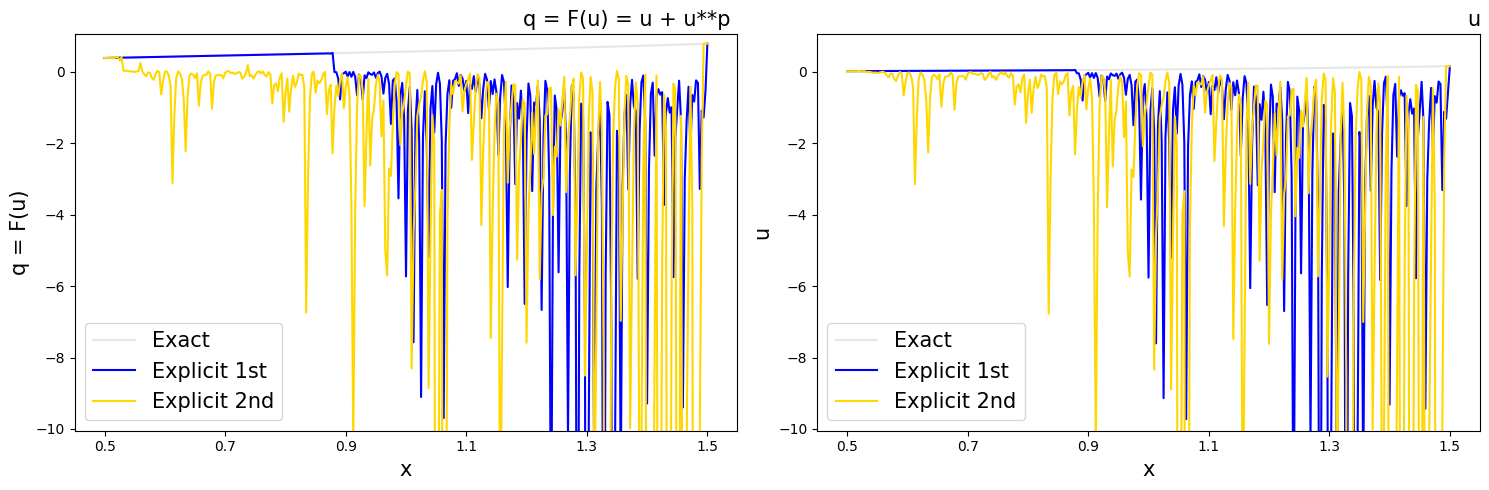}
    \caption{The visualization of numerical results of $q$ and $u$, obtained using the first (blue color) and second order (yellow color) accurate explicit numerical schemes, computed with a maximum Courant number $C_{max}=2$. The exact solution is displayed as a gray line. }
    \label{fig:ukazkavybuchu_cn4}
    \end{center}
\end{figure}

Additionally, the necessity to adhere to the stability restrictions that must be fulfilled when using explicit schemes is a concern; see Figure \ref{fig:ukazkavybuchu_cn4} for an illustration of blowup in numerical solution when using an explicit scheme with too large time step ($C_{max}=2$). Therefore, we will primarily focus on employing implicit numerical schemes from this point onward, which allow for much higher Courant numbers.

In the second table focused on the convergence order (Table \ref{TAB:Results_cn-10}), we provide the errors and the experimental order of convergence (EOC) for the higher Courant number exceeding the stability limits of explicit schemes, specifically we choose $N=M/20$ time steps, for $M=320, 640, 1280, 2560$, with the Courant number $C_{max} = 10$. This analysis focuses solely on the first order and second order implicit numerical schemes, highlighting their accuracy and behavior under no stability conditions. We provide results using 3 different choices of $p=1/4,1/2,3/4$ in (\ref{Freundlich}) and two different choices of $\omega$, first one constant value 0.5 and the second one using the WENO approximation \cite{zakova_numerical_2024}.

\begin{table}[h]
    \begin{center}
    \begin{tabular}{c c | c c | c c | c c }
    \hline
    \multicolumn{2}{c}{} & \multicolumn{6}{c}{First order}  \\
    \multicolumn{2}{c}{} & \multicolumn{2}{c}{$p=1/4$} & \multicolumn{2}{c}{$p=1/2$} & \multicolumn{2}{c}{$p=3/4$} \\
    \hline
    $M $ & $N$  & E                   & EOC  &  E                    & EOC  &  E                   & EOC  \\
    \hline
    320  & 16  & 2.21 $\cdot 10^{-3}$ & -    &  5.97 $\cdot 10^{-3}$ & -    & 1.68 $\cdot 10^{-2}$ & -  \\
    640  & 32  & 1.11 $\cdot 10^{-3}$ & 0.98 &  2.99 $\cdot 10^{-3}$ & 0.99 & 8.72 $\cdot 10^{-3}$ & 0.95    \\
    1280 & 64  & 5.61 $\cdot 10^{-4}$ & 0.99 &  1.50 $\cdot 10^{-3}$ & 0.99 & 4.34 $\cdot 10^{-3}$ & 1.00  \\
    2560 & 128 & 2.81 $\cdot 10^{-4}$ & 0.99 &  7.50 $\cdot 10^{-4}$ & 0.99 & 2.16 $\cdot 10^{-3}$ & 1.01   \\
    \hline 
    \hline
    \multicolumn{2}{c}{} & \multicolumn{6}{c}{Second order, $\omega=0.5$}  \\
    \multicolumn{2}{c}{}  & \multicolumn{2}{c}{$p=1/4$} & \multicolumn{2}{c}{$p=1/2$} & \multicolumn{2}{c}{$p=3/4$} \\
    \hline
    $M $ & $N$  & E                   & EOC  &  E                   & EOC  &  E                   & EOC  \\
    \hline
    320  & 16  & 4.02 $\cdot 10^{-5}$ & -    & 1.33 $\cdot 10^{-4}$ & -    & 1.44 $\cdot 10^{-3}$ & -     \\
    640  & 32  & 1.00 $\cdot 10^{-5}$ & 1.99 & 3.36 $\cdot 10^{-5}$ & 1.98 & 3.94 $\cdot 10^{-4}$ & 1.87 \\
    1280 & 64  & 2.52 $\cdot 10^{-6}$ & 1.99 & 8.22 $\cdot 10^{-6}$ & 2.03 & 6.42 $\cdot 10^{-5}$ & 2.61 \\
    2560 & 128 & 6.32 $\cdot 10^{-7}$ & 1.99 & 2.05 $\cdot 10^{-6}$ & 1.99 & 1.08 $\cdot 10^{-5}$ & 2.56 \\
    \hline 
    \hline
    \multicolumn{2}{c}{} & \multicolumn{6}{c}{High resolution WENO}  \\
    \multicolumn{2}{c}{}  & \multicolumn{2}{c}{$p=1/4$} & \multicolumn{2}{c}{$p=1/2$} & \multicolumn{2}{c}{$p=3/4$} \\
    \hline
    $M $ & $N$  & E                   & EOC  &  E                   & EOC  &  E                   & EOC  \\
    \hline
    320  & 16  & 4.61 $\cdot 10^{-5}$ & -    & 2.25 $\cdot 10^{-4}$ & -    & 2.24 $\cdot 10^{-3}$ & -     \\
    640  & 32  & 1.16 $\cdot 10^{-5}$ & 1.98 & 4.33 $\cdot 10^{-5}$ & 2.37 & 6.51 $\cdot 10^{-4}$ & 1.78 \\
    1280 & 64  & 2.93 $\cdot 10^{-6}$ & 1.99 & 9.52 $\cdot 10^{-6}$ & 2.18 & 1.55 $\cdot 10^{-5}$ & 2.06 \\
    2560 & 128 & 7.36 $\cdot 10^{-7}$ & 1.99 & 2.43 $\cdot 10^{-6}$ & 1.96 & 2.87 $\cdot 10^{-5}$ & 2.43 \\
    \hline 
    \end{tabular}
    \end{center}
    \caption{The numerical errors, and the EOCs using the first and second order accurate implicit schemes (with constant value of $\boldsymbol{\omega}$ and high resolution with WENO \cite{zakova_numerical_2024}), for a case with solutions without discontinuities, using higher maximum Courant number $C_{max} = 10$, with different choices of $p = 1/4,1/2,3/4$. }
    \label{TAB:Results_cn-10}
\end{table}

Secondly, we present the numerical results computed in the entire domain $x \in [0,5]$ for the initial condition (\ref{initialcondition}), which contains a discontinuity, and the exact solutions may be observed in Figures \ref{fig:Exact_pmensieako1}-\ref{fig:Exact_pvacsieako2}. As before, we compute the errors and the EOCs. However, when employing the second order accurate numerical scheme (\ref{scheme_2nd_implicit}), oscillations arise due to the presence of the discontinuity. To address this issue, WENO approximations were applied \cite{zakova_numerical_2024}, together with the high resolution modification of the scheme (\ref{scheme_2nd_implicit_withvel}). The results are summarized in Tables \ref{TAB:Results_cn-6_pmensieako1}-\ref{TAB:Results_cn-6_vacsieako2}, which provide a comparison for different choices of parameter $p$. 

\begin{figure}[h]
    \begin{center}
    \includegraphics[width=17cm]{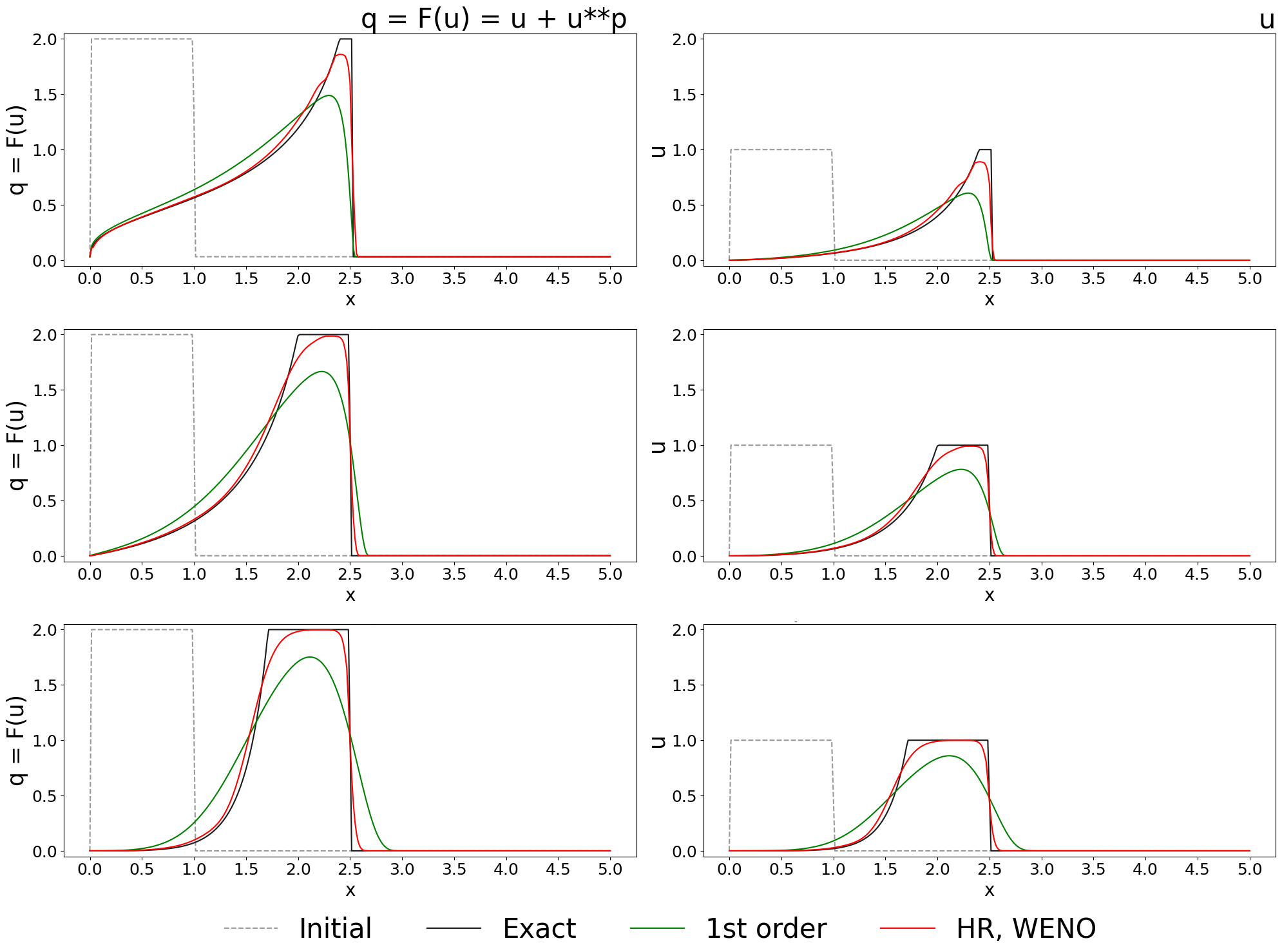}
    \caption{The initial condition, the exact solution, numerical solution - first and second order (high resolution with WENO approximations) - $q$ (left) and $u$ (right), for $p=1/4$ (top) $p=1/2$ (middle) $p=3/4$ (bottom) with the final time $T=3$, $M=320$, $C_{max}=6$.}
    \label{fig:numericlaWENO_pmensieako1}
    \end{center}
\end{figure}

\begin{figure}[h]
    \begin{center}
    \includegraphics[width=17cm]{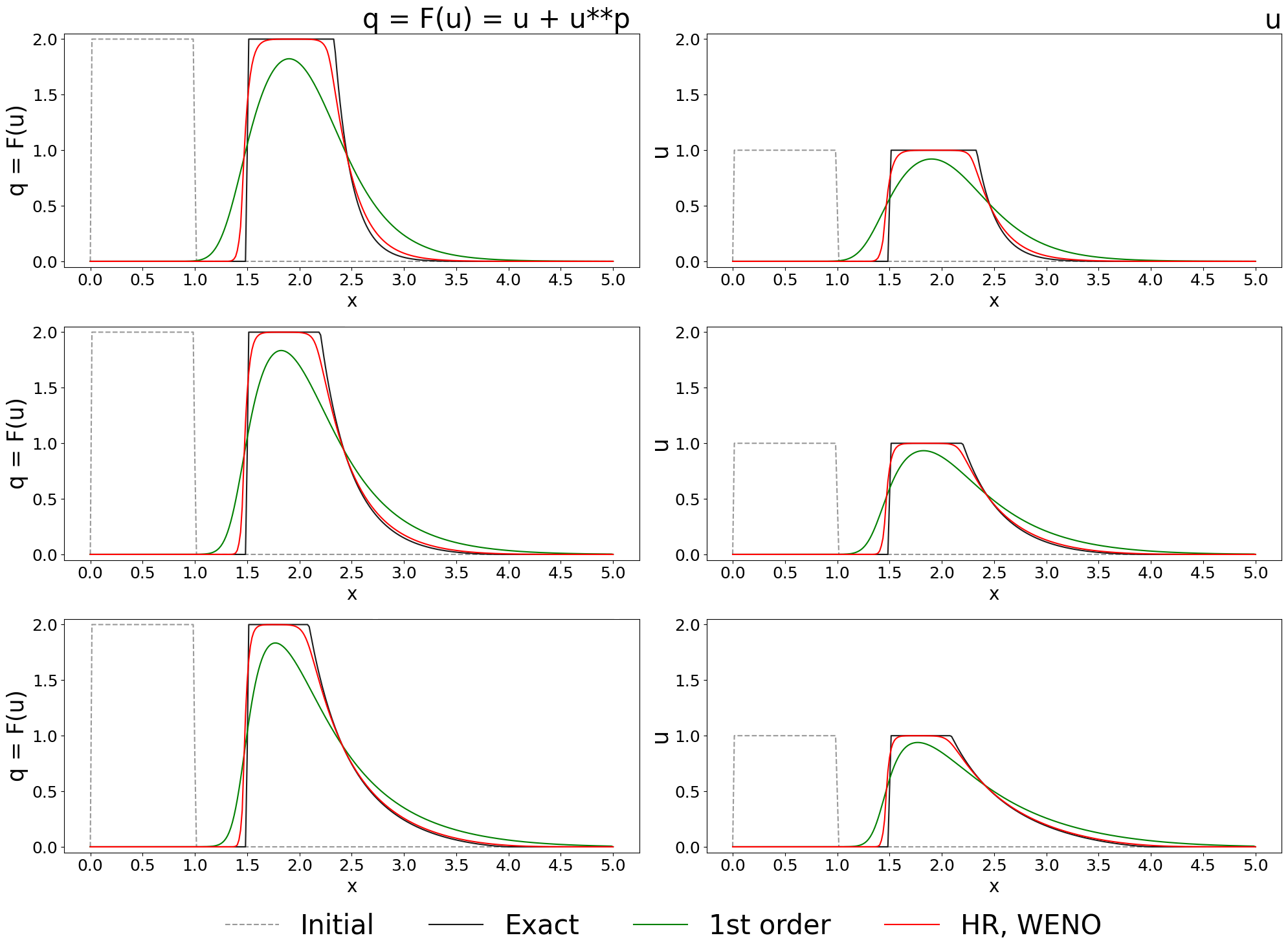}
    \caption{The initial condition, the exact solution, the numerical solution - first and second order (high resolution with WENO approximations) - $q$ (left) and $u$ (right), for $p=5/4$ (top) $p=3/2$ (middle) $p=7/4$ (bottom) with the final time $T=3$, $M=320$, $C_{max}=6$.}
    \label{fig:numericlaWENO_pvacsieako1}
    \end{center}
\end{figure}

\begin{figure}[h]
    \begin{center}
    \includegraphics[width=17cm]{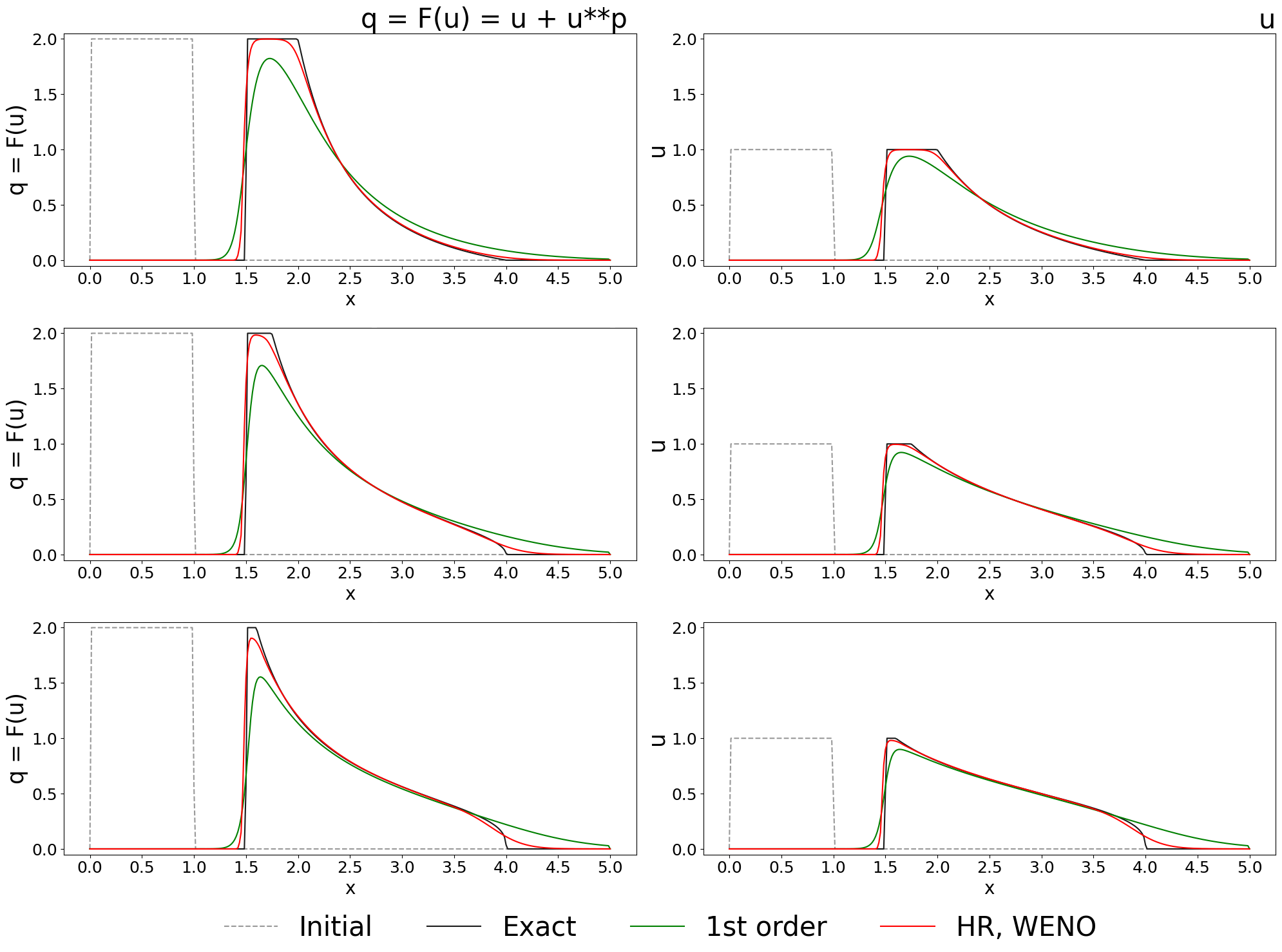}
    \caption{The initial condition, the exact solution, the numerical solution - first and second order (high resolution with WENO approximations) - $q$ (left) and $u$ (right), for $p=2$ (top) $p=3$ (middle) $p=4$ (bottom) with the final time $T=3$, $M=320$, $C_{max}=6$.}
    \label{fig:numericlaWENO_pvacsieako2}
    \end{center}
\end{figure}

We choose $N=M/10$ ($M=320, 640, 1280, 2560$) time steps with the Courant number being $C_{max}=4$. Visualizations of the initial condition, exact solution, and numerical results are provided for both $q$ and $u$, using the first order implicit and second order implicit (with WENO approximations) schemes. The results are displayed in Figure \ref{fig:numericlaWENO_pmensieako1} for $p=1/4,1/2,3/4$ (corresponding Table \ref{TAB:Results_cn-6_pmensieako1}), Figure \ref{fig:numericlaWENO_pvacsieako1} for $p=5/4,3/2,7/4$ (corresponding Table \ref{TAB:Results_cn-6_vacsieako1}), and in Figure \ref{fig:numericlaWENO_pvacsieako2} for $p=2,3,4$ (corresponding Table \ref{TAB:Results_cn-6_vacsieako2}).

\begin{table}[h]
    \begin{center}
    \begin{tabular}{c c | c c | c c | c c }
    \hline
    \multicolumn{2}{c}{} & \multicolumn{6}{c}{First order}  \\
    \multicolumn{2}{c}{} & \multicolumn{2}{c}{$p=1/4$} & \multicolumn{2}{c}{$p=1/2$} & \multicolumn{2}{c}{$p=3/4$} \\
    \hline
    $M $ & $N$  & E                   & EOC  &  E                    & EOC  &  E                   & EOC  \\
    \hline
    320  & 32  & 2.06 $\cdot 10^{-1}$ & -    & 2.71 $\cdot 10^{-1}$ & -    & 3.59 $\cdot 10^{-1}$ & -  \\
    640  & 64  & 1.45 $\cdot 10^{-1}$ & 0.50 & 1.76 $\cdot 10^{-1}$ & 0.62 & 2.32 $\cdot 10^{-1}$ & 0.63    \\
    1280 & 128 & 9.56 $\cdot 10^{-2}$ & 0.59 & 1.10 $\cdot 10^{-1}$ & 0.67 & 1.44 $\cdot 10^{-1}$ & 0.68  \\
    2560 & 256 & 6.33 $\cdot 10^{-2}$ & 0.60 & 6.75 $\cdot 10^{-2}$ & 0.70 & 8.82 $\cdot 10^{-2}$ & 0.71   \\
    \hline 
    \hline
    \multicolumn{2}{c}{} & \multicolumn{6}{c}{High resolution WENO}  \\
    \multicolumn{2}{c}{}  & \multicolumn{2}{c}{$p=1/4$} & \multicolumn{2}{c}{$p=1/2$} & \multicolumn{2}{c}{$p=3/4$} \\
    \hline
    $M $ & $N$  & E                   & EOC  &  E                   & EOC  &  E                   & EOC   \\
    \hline
    320  & 32  & 6.94 $\cdot 10^{-2}$ & -    & 7.81 $\cdot 10^{-2}$ & -    & 9.25 $\cdot 10^{-2}$ & -     \\
    640  & 64  & 4.06 $\cdot 10^{-2}$ & 0.77 & 4.03 $\cdot 10^{-2}$ & 0.95 & 4.83 $\cdot 10^{-2}$ & 0.94  \\
    1280 & 128 & 2.14 $\cdot 10^{-2}$ & 0.92 & 2.06 $\cdot 10^{-2}$ & 0.97 & 2.50 $\cdot 10^{-2}$ & 0.95  \\
    2560 & 256 & 1.09 $\cdot 10^{-2}$ & 0.97 & 1.04 $\cdot 10^{-2}$ & 0.99 & 1.27 $\cdot 10^{-2}$ & 0.97  \\
    \hline 
    \end{tabular}
    \end{center}
    \caption{The numerical errors, and the EOCs using the first and second order (high resolution with WENO approximations) accurate implicit schemes, for a case with a discontinuous initial condition (\ref{InitialCondition_discontinuous}), using the maximum Courant number $C_{max} = 6$, with different choices of $p = 1/4,1/2,3/4$. }
    \label{TAB:Results_cn-6_pmensieako1}
\end{table}

\begin{table}[h]
    \begin{center}
    \begin{tabular}{c c | c c | c c | c c }
    \hline
    \multicolumn{2}{c}{} & \multicolumn{6}{c}{First order}  \\
    \multicolumn{2}{c}{} & \multicolumn{2}{c}{$p=5/4$} & \multicolumn{2}{c}{$p=3/2$} & \multicolumn{2}{c}{$p=7/4$} \\
    \hline
    $M $ & $N$  & E                   & EOC  &  E                    & EOC  &  E                   & EOC  \\
    \hline
    320  & 32  & 3.94 $\cdot 10^{-1}$ & -    & 3.34 $\cdot 10^{-1}$ & -    & 2.94 $\cdot 10^{-1}$ & -     \\
    640  & 64  & 2.53 $\cdot 10^{-1}$ & 0.63 & 2.03 $\cdot 10^{-1}$ & 0.71 & 1.74 $\cdot 10^{-1}$ & 0.75  \\
    1280 & 128 & 1.58 $\cdot 10^{-1}$ & 0.68 & 1.20 $\cdot 10^{-1}$ & 0.75 & 1.01 $\cdot 10^{-1}$ & 0.78  \\
    2560 & 256 & 9.55 $\cdot 10^{-2}$ & 0.72 & 6.99 $\cdot 10^{-2}$ & 0.77 & 5.83 $\cdot 10^{-2}$ & 0.79  \\
    \hline 
    \hline
    \multicolumn{2}{c}{} & \multicolumn{6}{c}{High resolution WENO}  \\
    \multicolumn{2}{c}{} & \multicolumn{2}{c}{$p=5/4$} & \multicolumn{2}{c}{$p=3/2$} & \multicolumn{2}{c}{$p=7/4$} \\
    \hline
    $M $ & $N$  & E                   & EOC  &  E                   & EOC  &  E                   & EOC   \\
    \hline
    320  & 32  & 1.08 $\cdot 10^{-1}$ & -    & 9.12 $\cdot 10^{-2}$ & -    & 8.29 $\cdot 10^{-2}$ & -     \\
    640  & 64  & 5.54 $\cdot 10^{-2}$ & 0.96 & 4.59 $\cdot 10^{-2}$ & 0.99 & 4.15 $\cdot 10^{-2}$ & 0.99  \\
    1280 & 128 & 2.81 $\cdot 10^{-2}$ & 0.98 & 2.30 $\cdot 10^{-2}$ & 0.99 & 2.08 $\cdot 10^{-2}$ & 0.99  \\
    2560 & 256 & 1.41 $\cdot 10^{-2}$ & 0.99 & 1.15 $\cdot 10^{-2}$ & 0.99 & 1.04 $\cdot 10^{-2}$ & 0.99  \\
    \hline 
    \end{tabular}
    \end{center}
    \caption{The numerical errors, and the EOCs using the first and second order (high resolution with WENO approximations) accurate implicit schemes, for a case with a discontinuous initial condition (\ref{InitialCondition_discontinuous}), using the maximum Courant number $C_{max} = 6$, with different choices of $p = 5/4,3/2,7/4$. }
    \label{TAB:Results_cn-6_vacsieako1}
\end{table}

\begin{table}[h]
    \begin{center}
    \begin{tabular}{c c | c c | c c | c c }
    \hline
    \multicolumn{2}{c}{} & \multicolumn{6}{c}{First order}  \\
    \multicolumn{2}{c}{} & \multicolumn{2}{c}{$p=2$} & \multicolumn{2}{c}{$p=3$} & \multicolumn{2}{c}{$p=4$} \\
    \hline
    $M $ & $N$  & E                   & EOC  &  E                    & EOC  &  E                   & EOC  \\
    \hline
    320  & 32  & 2.66 $\cdot 10^{-1}$ & -    & 2.06 $\cdot 10^{-1}$ & -    & 2.03 $\cdot 10^{-1}$ & -     \\
    640  & 64  & 1.56 $\cdot 10^{-1}$ & 0.76 & 1.21 $\cdot 10^{-1}$ & 0.77 & 1.27 $\cdot 10^{-1}$ & 0.67  \\
    1280 & 128 & 9.03 $\cdot 10^{-2}$ & 0.79 & 6.88 $\cdot 10^{-2}$ & 0.81 & 7.87 $\cdot 10^{-2}$ & 0.69  \\
    2560 & 256 & 5.15 $\cdot 10^{-2}$ & 0.80 & 3.84 $\cdot 10^{-2}$ & 0.83 & 4.89 $\cdot 10^{-2}$ & 0.68  \\
    \hline 
    \hline
    \multicolumn{2}{c}{} & \multicolumn{6}{c}{High resolution WENO}  \\
    \multicolumn{2}{c}{} & \multicolumn{2}{c}{$p=2$} & \multicolumn{2}{c}{$p=3$} & \multicolumn{2}{c}{$p=4$} \\
    \hline
    $M $ & $N$  & E                   & EOC  &  E                   & EOC  &  E                   & EOC   \\
    \hline
    320  & 32  & 7.81 $\cdot 10^{-2}$ & -    & 6.02 $\cdot 10^{-2}$ & -    & 7.27 $\cdot 10^{-2}$ & -     \\
    640  & 64  & 3.91 $\cdot 10^{-2}$ & 0.99 & 3.02 $\cdot 10^{-2}$ & 0.99 & 3.81 $\cdot 10^{-2}$ & 0.93  \\
    1280 & 128 & 1.95 $\cdot 10^{-2}$ & 0.99 & 1.51 $\cdot 10^{-2}$ & 0.99 & 1.99 $\cdot 10^{-2}$ & 0.94  \\
    2560 & 256 & 9.78 $\cdot 10^{-3}$ & 0.99 & 7.59 $\cdot 10^{-3}$ & 0.99 & 1.03 $\cdot 10^{-2}$ & 0.94  \\
    \hline 
    \end{tabular}
    \end{center}
    \caption{The numerical errors, and the EOCs using the first and second order (high resolution with WENO approximations) accurate implicit schemes, for a case with a discontinuous initial condition (\ref{InitialCondition_discontinuous}), using the maximum Courant number $C_{max} = 6$, with different choices of $p = 2, 3, 4$.}
    \label{TAB:Results_cn-6_vacsieako2}
\end{table}

\clearpage

\subsection{One dimensional experiment with variable velocity}
\label{subsec_1D_withvel}

\begin{figure}[h!]
    \begin{center}
    \includegraphics[width=17cm, height=14cm]{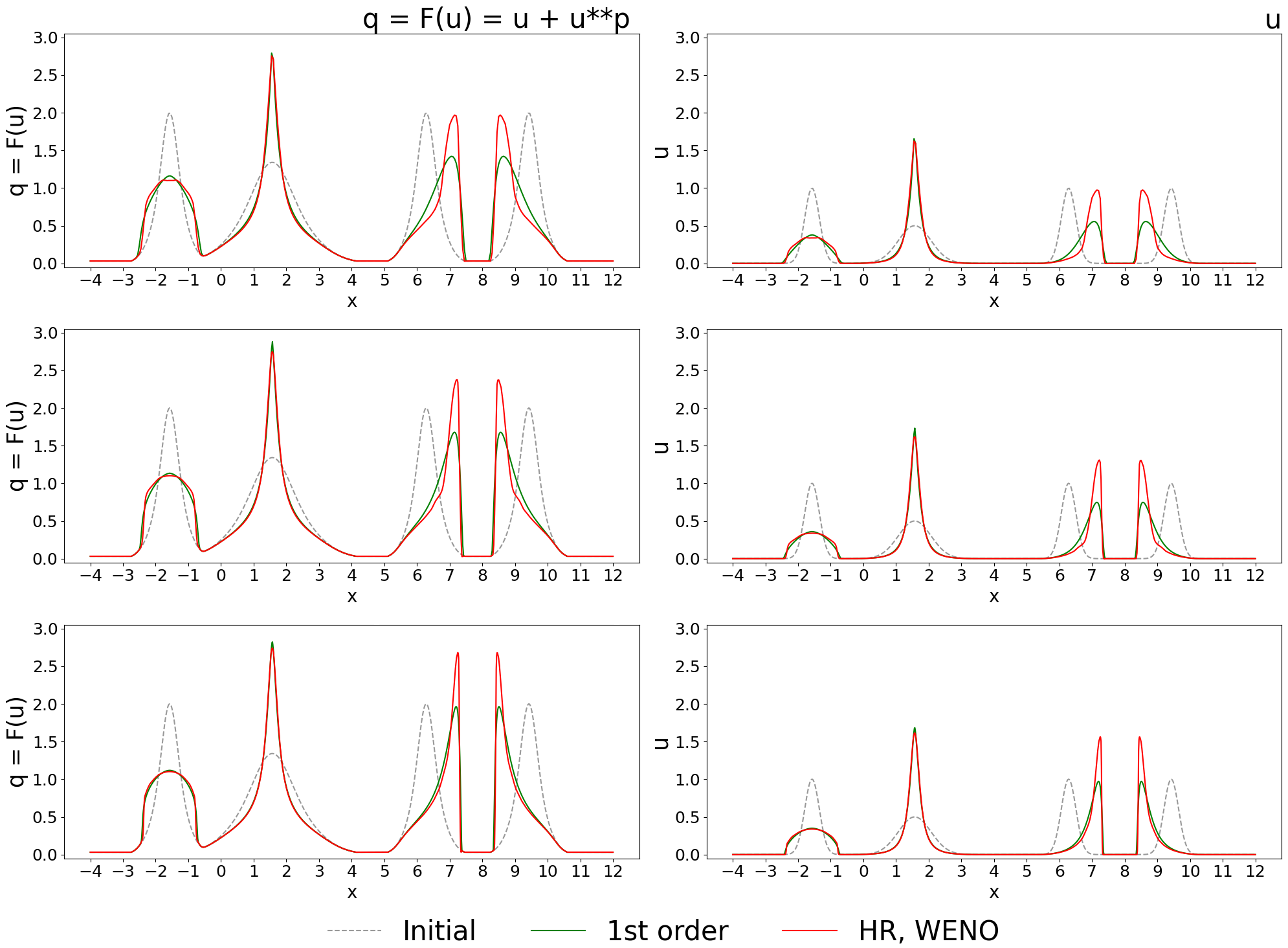}
    \caption{The comparison of the numerical solution ($q$ on the left side and $u$ on the right) first order (green color), high resolution second order with  WENO approximation (red color), for $p=1/4$ for different choices of $M = 320$ (top) $M = 640$ (middle) $M = 1280$ (bottom) for the final time $T=1.5$. The initial condition is shown as a gray dashed line.}
    \label{fig:pexp0.25_results_cn7.5_WENO_320az1280}
    \end{center}
\end{figure}

\begin{figure}[h!]
    \begin{center}
    \includegraphics[width=17cm, height = 14cm]{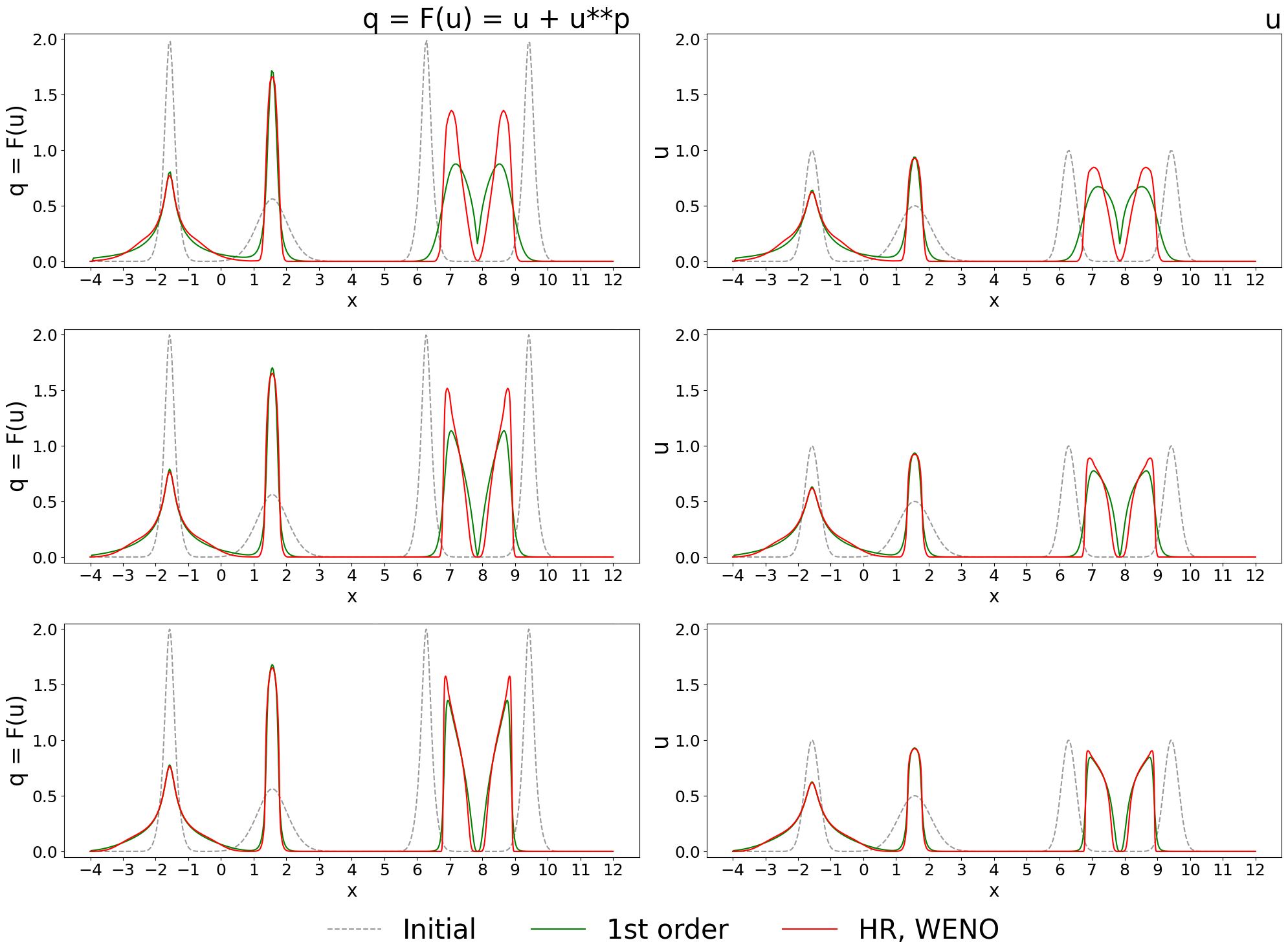}
    \caption{The comparison of the numerical solution ($q$ on the left side and $u$ on the right) first order (green color), high resolution second order with  WENO approximation (red color), for $p=4$ for different choices of $M = 320$ (top) $M = 640$ (middle) $M = 1280$ (bottom) for the final time $T=1.5$. The initial condition is shown as a gray dashed line.}
    \label{fig:pexp4_results_cn7.5_WENO_320az1280}
    \end{center}
\end{figure}

The next experiment in one dimension demonstrates a visual comparison (Figures \ref{fig:pexp0.25_results_cn7.5_WENO_320az1280}, \ref{fig:pexp4_results_cn7.5_WENO_320az1280}) of the results using the variable velocity $v(x)=\cos(x)$, which can take on both positive and negative values, as described in equation (\ref{eq_general_withvel}); therefore, it is necessary to perform at least two sweeps (\ref{sweepiter1})-(\ref{sweepiter2}) of the fast sweeping method. In this case, we choose $x \in [-4,11]$ and $T=1.5$ with $N=M/80$ ($M$ defined later) and $C_{max} = 7.5$ (\ref{Courantnumber}). 

The initial condition is represented by a sum of four Gaussian functions as
\begin{equation}
    \label{4gausians_1D}
    u^0(x) = e^{-10(x + \frac{\pi}{2})^2 } + 0.5 e^{-2(x - \frac{\pi}{2})^2} + e^{-10(x - 2\pi)^2} + e^{-10(x - 3\pi)^2 }  \,
\end{equation}
with the boundary condition $u(-4,t)=0, u(12,t)=0 $. To utilize the full stencil of the scheme for every point in the interior grid, we extend the zero boundary condition to the neighboring points of the boundaries.

Since we do not know the exact solution for this case, we present only the numerical solution in the final time $T$ using three different domain meshes, which are sequentially refined, namely $M=320 \, (N=4), 640 \,(N=8), 1280\, (N=16)$. The visual results can be found in Figure \ref{fig:pexp0.25_results_cn7.5_WENO_320az1280} for $p=1/4$ and Figure \ref{fig:pexp4_results_cn7.5_WENO_320az1280} for $p=4$. Notice the notable differences and improvements when the mesh is refined and when the second order accurate scheme is employed.

\clearpage

\subsection{Two dimensional experiment with velocity field}
\label{subsec_2D_withvel}

The last numerical experiment will focus on the two dimensional case, with the initial condition, which again consists of a sum of four Gaussian functions,  defined as 
\begin{align}
    \label{2Dgaussians}
    u^0(x,y) = \,\,  & e^{-50\left((x-0.5)^2 + (y+0.5)^2\right)} + e^{-50\left((x-0.5)^2 + (y-0.5)^2\right)}  \nonumber \\
              + \,\, & e^{-50\left((x+0.5)^2 + (y+0.5)^2\right)} + e^{-50\left((x+0.5)^2 + (y-0.5)^2\right)}
\end{align}
for $x,y\in[-1, 1]$ and we use the final time $T=1/4$ with $N=M/10$ time steps ($M$ defined later). With the boundary conditions set to zero at every boundary of the computational domain, we also extend the zero boundary condition to neighboring points outside the boundaries to ensure that the full stencil of the scheme can be applied to every point in the interior grid.

\begin{wrapfigure}{r}{0.3\textwidth}
    \centering
    \vspace{-0.5cm}
    \includegraphics[width=0.9\linewidth]{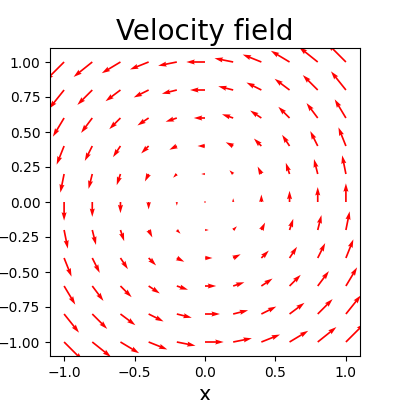}
    \vspace{-0.4cm}
    \caption{Velocity field $\Vec{v} = (-2\pi y,2\pi x)$. }
    \vspace{0.7cm}
    \label{FIG:velocityfield}
\end{wrapfigure} 

We will let the initial condition rotate, using the equation (\ref{2D_withvel}) with the divergence free velocity field $\Vec{v} = (-2\pi y,2\pi x)$ (see Figure \ref{FIG:velocityfield}), which can be either positive or negative depending on the specific location, therefore, it is necessary to perform all four sweeps of the fast sweeping method (\ref{2D_sweeps}). In this case, since the exact solution is again unknown, we will provide only a visual demonstration of the progression of the numerical solution as the mesh is refined, specifically $M=80 \,(N=8), 160\, (N=16), 320\, (N=32)$, sequentially. 

\begin{figure}[h!]
    \begin{center}
    \includegraphics[width=17cm]{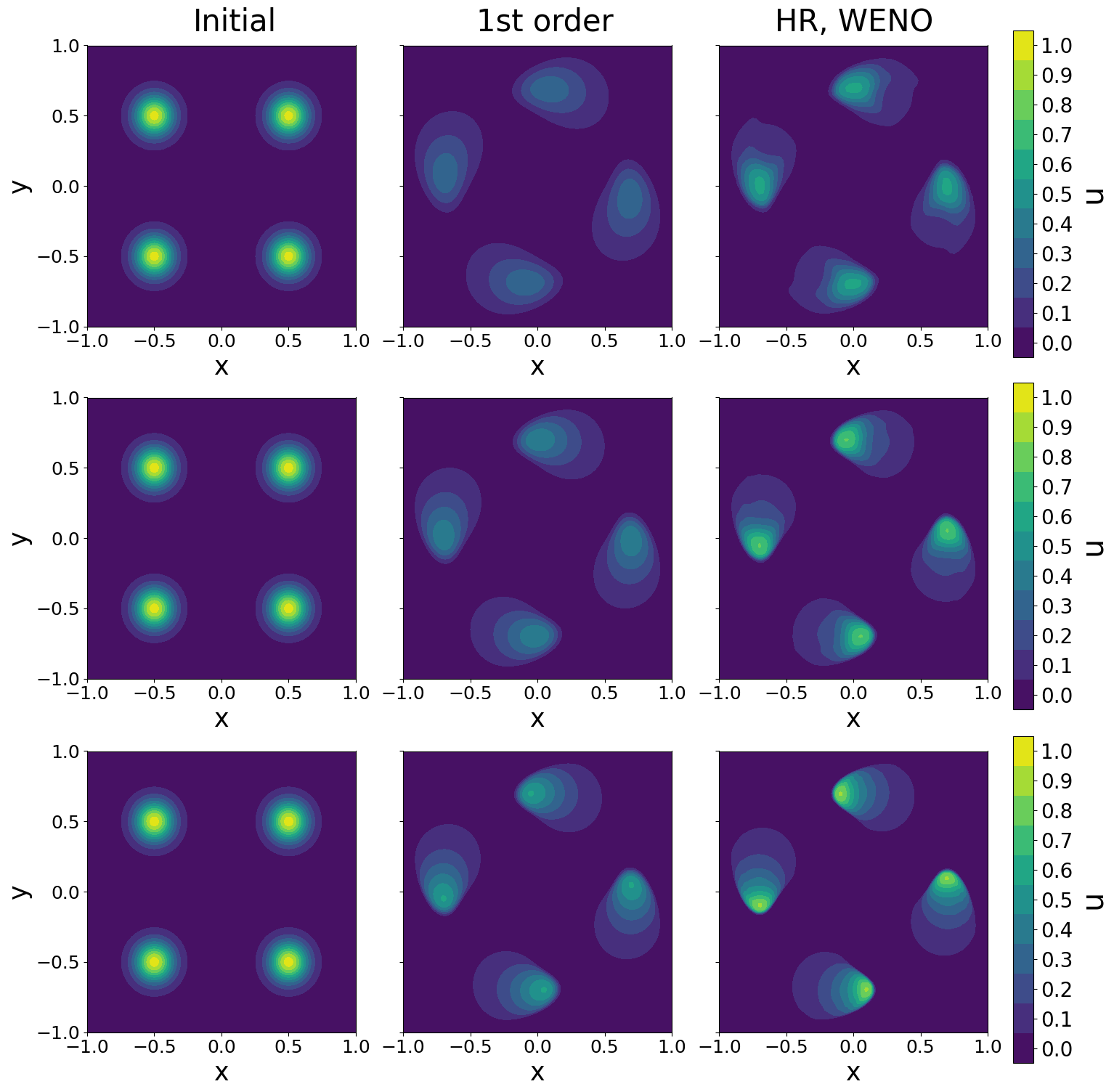}
    \caption{The comparison of the contours of the numerical solution $u$ - first order (middle), high resolution second order with  WENO approximation (right), for $p=1/2$ for different choices of $M = 80$ (top) $M = 160$ (middle) $M = 320$ (bottom) for the final time $T=1/4$. The initial condition is shown on the left.}
    \label{fig:2D_4gausiany_pexp0.5_Cn7.85_FINAL_80az320}
    \end{center}
\end{figure}

\begin{figure}[h!]
    \begin{center}
    \includegraphics[width=17cm]{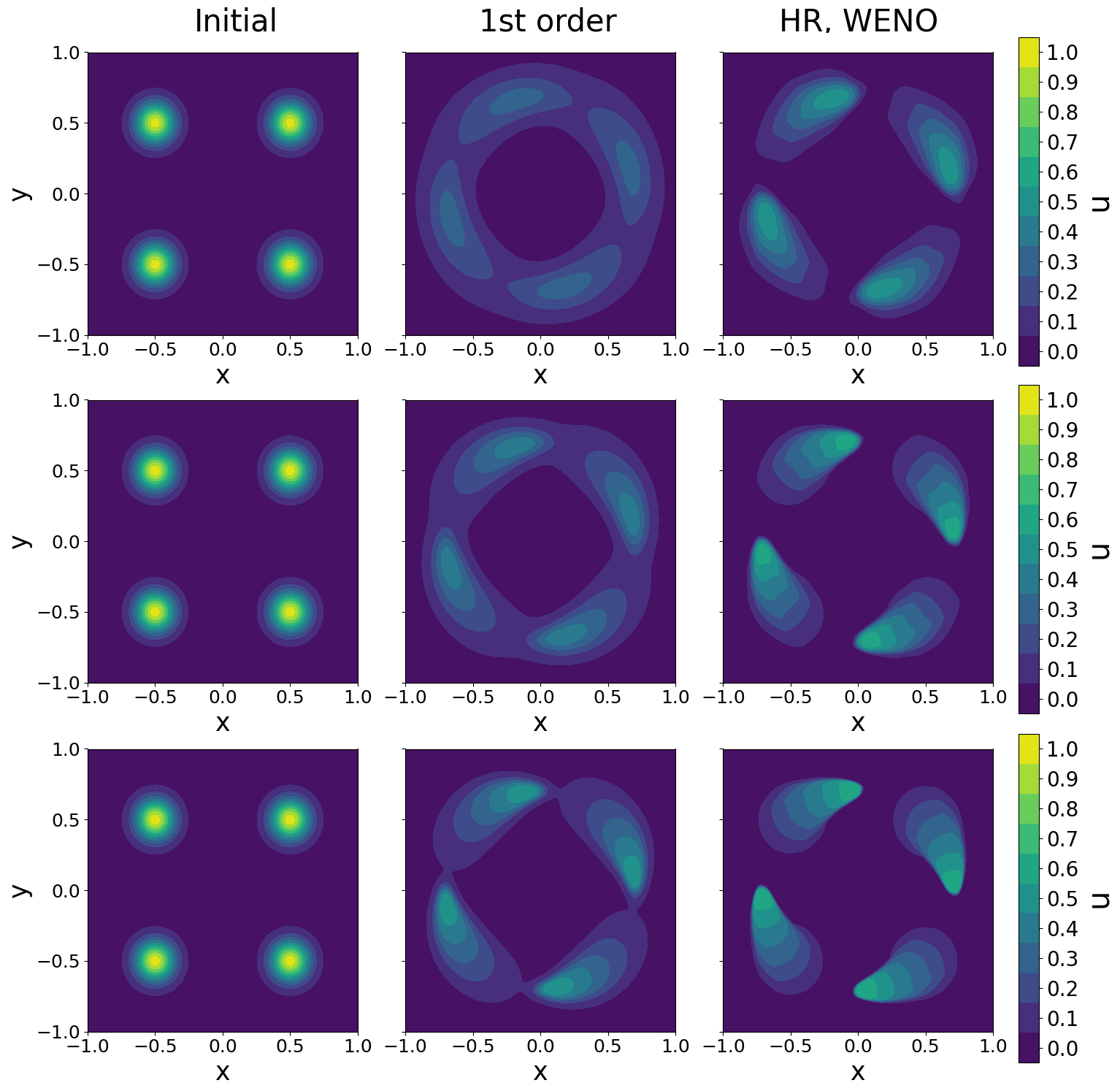}
    \caption{The comparison of the contours of the numerical solution $u$ - first order (middle), high resolution second order with  WENO approximation (right), for $p=3$ for different choices of $M = 80$ (top) $M = 160$ (middle) $M = 320$ (bottom) for the final time $T=1/4$. The initial condition is shown on the left.}
    \label{fig:2D_4gausiany_pexp3_Cn7.85_FINAL_80az320}
    \end{center}
\end{figure}

The contours of the numerical results can be found in Figure \ref{fig:2D_4gausiany_pexp0.5_Cn7.85_FINAL_80az320} for $p=1/2$ and in Figure \ref{fig:2D_4gausiany_pexp3_Cn7.85_FINAL_80az320} for $p=3$, which compare the visual results using the first order accurate and second order accurate high resolution (using WENO approximations) schemes. The figures demonstrate notable differences and improvements when the second order accurate scheme is employed and when the mesh is refined.

For the two dimensional case, we define two maximum Courant numbers, one for each of the directions,  
\begin{align}
    \label{Courantnumber_2D}
    &C_{max}^x = 
    \frac{\tau}{h} \max \limits_{i,j} \Bigl\{ \max \bigl\{ |v(x_{i-1/2},y_j)|, |v(x_{i+1/2},y_j)| \bigr\}  \Bigr\} \,, \nonumber \\
    &C_{max}^y =  
    \frac{\tau}{h} \max \limits_{i,j} \Bigl\{ \max \bigl\{ |w(x_{i},y_{j-1/2})|, |w(x_{i},y_{j+1/2})| \bigr\}  \Bigr\} \,.
\end{align}

The maximum Courant number in this case is the same for each of the directions $C_{max}^x=C_{max}^y=7.85$.

\clearpage
\section{Conclusion}
\label{sec-conclude}

This study presents implicit numerical methods for solving nonlinear transport problems characterized by various sorption isotherms, particularly focusing on the Freundlich isotherm. Through the development of both first and second order accurate numerical schemes, including implicit and explicit methods, we have demonstrated the effectiveness of these approaches in accurately modeling the transport phenomena under different conditions.

The results from our numerical experiments highlight the advantages of implicit schemes over explicit ones, particularly in terms of stability and computational time efficiency when dealing with higher Courant numbers. The implementation of high resolution methods has proven to be important in minimizing oscillations that arise near discontinuities, consequently increasing the overall plausibility of the numerical solutions. The numerical convergence rates observed suggest that both first and second order schemes can achieve satisfactory accuracy, with the second order methods providing significant improvements in precision.

\rv{
Overall, this work contributes to the theoretical understanding of nonlinear transport phenomena by developing and analyzing compact implicit high resolution numerical schemes that effectively address the nonlinearities present in transport equations with sorption isotherms. In addition, it provides practical numerical tools applicable to various fields, such as environmental hydrology and materials science. The paper demonstrates that these schemes achieve higher accuracy without the restrictive time step limitations of explicit methods and that the use of limiters in high resolution schemes ensures physically realistic solutions even in the presence of discontinuities.} To offer a significant practical value of the numerical schemes presented, future research must explore their development for a system of coupled transport equations \cite{donat_weno_2024, schmidt-traubPreparativeChromatography2020, frolkovicNumericalSimulationContaminant2016}.

\printbibliography

@article{frolkovic_semi-analytical_2006,
	title = {Semi-analytical solutions of a contaminant transport equation with nonlinear sorption in {1D}},
	volume = {10},
	issn = {1420-0597, 1573-1499},
        language = {en},
	number = {3},
	journal = {Computational Geosciences},
	author = {Frolkovič, Peter and Kačur, Jozef},
	month = sep,
	year = {2006},
	pages = {279--290},
}

@article{donatImplicitExplicitWENO2018,
  title = {Implicit--{{Explicit WENO}} Scheme for the Equilibrium Dispersive Model of Chromatography},
  author = {Donat, R. and Guerrero, F. and Mulet, P.},
  year = {2018},
  journal = {Applied Numerical Mathematics},
  volume = {123},
  pages = {22--42},
}

@article{burger_linearly_2018,
	title = {Linearly implicit-explicit schemes for the equilibrium dispersive model of chromatography},
	volume = {317},
	issn = {0096-3003},
	journal = {Applied Mathematics and Computation},
	author = {Bürger, Raimund and Mulet, Pep and Rubio, Lihki and Sepúlveda, Mauricio},
	month = jan,
	year = {2018},
	pages = {172--186},
}

@article{donat_weno_2024,
	title = {{WENO} scheme on characteristics for the equilibrium dispersive model of chromatography with generalized {Langmuir} isotherms},
	volume = {201},
	issn = {0168-9274},
	journal = {Applied Numerical Mathematics},
	author = {Donat, R. and Martí, M. C. and Mulet, P.},
	month = jul,
	year = {2024},
	pages = {247--264},
}

@book{leveque_finite_2004,
	edition = {2nd},
	title = {Finite {Volume} {Methods} for {Hyperbolic} {Problems}},
	language = {en},
	publisher = {Cambridge UP},
	author = {Leveque, Randall J},
	year = {2004},

}

@article{frolkovic2023high,
  title = {High Resolution Compact Implicit Numerical Scheme for Conservation Laws},
  author = {Frolkovi{\v c}, Peter and {\v Z}erav{\'y}, Michal},
  year = {2023},
  journal = {Applied Mathematics and Computation},
  volume = {442},
  pages = {127720}
}

@misc{zakova_numerical_2024,
	title = {Numerical solution of two dimensional scalar conservation laws using compact implicit numerical schemes},
	language = {en},
	publisher = {arXiv},
	author = {Žáková, Dagmar and Frolkovič, Peter},
	month = aug,
	year = {2024},
	note = {arXiv:2407.05275 [math]},	
}

@Techreport{CS-R9526,
    title= {Python tutorial},    
    author = {G. van Rossum}, 
    institution= {Centrum voor Wiskunde en Informatica (CWI)},    
    year= {1995},    
    month={5}  
}

@incollection{shu_essentially_1998,
	address = {Berlin, Heidelberg},
	series = {Lecture {Notes} in {Mathematics}},
	title = {Essentially non-oscillatory and weighted essentially non-oscillatory schemes for hyperbolic conservation laws},
	booktitle = {Advanced {Numerical} {Approximation} of {Nonlinear} {Hyperbolic} {Equations}},
	publisher = {Springer},
	author = {Shu, Chi-Wang},
	year = {1998},
}

@article{qiu_finite_2003,
	title = {Finite {Difference} {WENO} {Schemes} with {Lax}--{Wendroff}-{Type} {Time} {Discretizations}},
	volume = {24},
  number={6},
  pages={2185--2198},
	journal = {SIAM J. Sci. Comp.},
	author = {Qiu, Jianxian and Shu, Chi-Wang},
	year = {2003},
}

@incollection{AACEE_2023_Zakova,
  author      = "Dagmar Zakova",
  title       = "Second order accurate compact implicit numerical scheme for solving linear advection equation in two dimensions",
  editor      = "Bisták, Andrej",
  booktitle   = "Advances in architectural, civil and environmental engineering 33rd annual {PhD} student conference on applied mathematics, building technology, geodesy and cartography, landscaping, theory and environmental technology of buildings, theory and structures of buildings, theory and structures of civil engineering works, water resources engineering: {October} 25th 2023, {Bratislava}, {Slovakia}",
  publisher   = "SPEKTRUM STU",
  address     = "Bratislava",
  year        = 2023,
  pages       = "40-47",
}

@article{duraisamy_implicit_2007,
	title = {Implicit {Scheme} for {Hyperbolic} {Conservation} {Laws} {Using} {Nonoscillatory} {Reconstruction} in {Space} and {Time}},
	volume = {29},
	number = {6},
	journal = {SIAM J. Sci. Comput.},
	author = {Duraisamy, Karthikeyan and Baeder, James D.},
	month = jan,
	year = {2007},
	
}

@article{javeedEfficientAccurateNumerical2011,
  title = {Efficient and Accurate Numerical Simulation of Nonlinear Chromatographic Processes},
  author = {Javeed, Shumaila and Qamar, Shamsul and {Seidel-Morgenstern}, Andreas and Warnecke, Gerald},
  year = {2011},
  journal = {Computers \& Chemical Engineering},
  volume = {35},
  number = {11},
  pages = {2294--2305},
}

@article{qamarAnalyticalSolutionsMoment2014,
  title = {Analytical Solutions and Moment Analysis of General Rate Model for Linear Liquid Chromatography},
  author = {Qamar, Shamsul and Nawaz Abbasi, Javeria and Javeed, Shumaila and {Seidel-Morgenstern}, Andreas},
  year = {2014},
  journal = {Chemical Engineering Science},
  volume = {107},
  pages = {192--205},
}

@article{zafarDiscontinuousGalerkinScheme2021,
  title = {Discontinuous {{Galerkin}} Scheme for Solving Non-Isothermal and Non-Equilibrium Model of Liquid Chromatography},
  author = {Zafar, Shireen and Perveen, Sadia and Qamar, Shamsul},
  year = {2021},
  journal = {Journal of Liquid Chromatography \& Related Technologies},
  volume = {44},
  number = {1-2},
  pages = {52--69},
}

@article{kaczmarskiAnalyticalNumericalSolutions2024,
  title = {Analytical and Numerical Solutions of Linear and Nonlinear Chromatography Column Models},
  author = {Kaczmarski, Krzysztof and Szukiewicz, Miros{\l}aw Krzysztof},
  year = {2024},
  journal = {Acta Chromatographica},
  volume = {-1},
  number = {aop},
}

@article{kacurSolutionContaminantTransport2005,
  title = {Solution of Contaminant Transport with Equilibrium and Non-Equilibrium Adsorption},
  author = {Ka{\v c}ur, J. and Malengier, B. and Reme{\v s}{\'i}kov{\'a}, M.},
  year = {2005},
  journal = {Computer Methods in Applied Mechanics and Engineering},
  series = {Selected Papers from the 11th {{Conference}} on {{The Mathematics}} of {{Finite Elements}} and {{Applications}}},
  volume = {194},
  number = {2},
  pages = {479--489},
}

@article{frolkovicNumericalSimulationContaminant2016,
  title = {Numerical Simulation of Contaminant Transport in Groundwater Using Software Tools of r3t},
  author = {Frolkovi{\v c}, Peter and Lampe, Michael and Wittum, Gabriel},
  year = {2016},
  journal = {Computing and Visualization in Science},
  volume = {18},
  number = {1},
  pages = {17--29},
}

@book{guMathematicalModelingScaleUp2015,
  title = {Mathematical {{Modeling}} and {{Scale-Up}} of {{Liquid Chromatography}}: {{With Application Examples}}},
  shorttitle = {Mathematical {{Modeling}} and {{Scale-Up}} of {{Liquid Chromatography}}},
  author = {Gu, Tingyue},
  year = {2015},
  publisher = {Springer},
  googlebooks = {PiryBwAAQBAJ},
  isbn = {978-3-319-16145-7}
}

@book{schmidt-traubPreparativeChromatography2020,
  title = {Preparative {{Chromatography}}},
  author = {{Schmidt-Traub}, H. and Schulte, Michael and {Seidel-Morgenstern}, Andreas},
  year = {2020},
  publisher = {John Wiley \& Sons},
  googlebooks = {7H\_PDwAAQBAJ},
  isbn = {978-3-527-34486-4}
}

@article{brusseau1995effect,
  title={The effect of nonlinear sorption on transformation of contaminants during transport in porous media},
  author={Brusseau, Mark L},
  journal={Journal of Contaminant Hydrology},
  volume={17},
  number={4},
  pages={277--291},
  year={1995},
  publisher={Elsevier}
}

@article{zhao2005fast,
  title={A fast sweeping method for eikonal equations},
  author={Zhao, Hongkai},
  journal={Mathematics of computation},
  volume={74},
  number={250},
  pages={603--627},
  year={2005}
}
\end{document}